\documentclass[]{article}
\usepackage{graphicx}
\usepackage{epstopdf}
\usepackage{amsfonts}
\usepackage{amsmath}
\usepackage{amssymb}
\usepackage{fancyhdr}
\usepackage{titlesec}
\usepackage{indentfirst}
\usepackage{booktabs}
\usepackage{verbatim}
\usepackage{color}
\usepackage{amsthm}
\usepackage{subfigure}
\usepackage{pdflscape}
\usepackage{todonotes}
\usepackage{bm}
\usepackage{amsmath,amssymb,amscd,mathrsfs,amsthm}
\usepackage{mathrsfs}
\usepackage{float}
\usepackage{url}
\usepackage{indentfirst}
\usepackage[colorlinks,linkcolor=black,anchorcolor=black,citecolor=black]{hyperref}

\usepackage[page,header]{appendix}
\usepackage{titletoc}
\usepackage{cleveref}

\newcommand{\rd}{\,\mathrm{d}}

\newcommand{\p}{\partial}

\newcommand{\fl}[2]{\frac{#1}{#2}}

\newcommand{\dt}{\Delta}

\newcommand{\be}{\begin{equation}}
\newcommand{\ee}{\end{equation}}
\newcommand{\bal}{\begin{aligned}}
\newcommand{\eal}{\end{aligned}}
\newcommand{\ba}{\begin{array}}
\newcommand{\ea}{\end{array}}
\newcommand{\bea}{\begin{eqnarray}}
\newcommand{\eea}{\end{eqnarray}}
\newcommand{\beas}{\begin{eqnarray*}}
\newcommand{\eeas}{\end{eqnarray*}}
\newcommand{\beit}{\begin{itemize}}
\newcommand{\eit}{\end{itemize}}

\newcommand{\bmat}[1]{\begin{bmatrix} #1 \end{bmatrix}}
\newcommand{\mat}[1]{\boldsymbol{#1}}
\renewcommand{\vec}[1]{\boldsymbol{\mathrm{#1}}}

\providecommand{\mK}{\ensuremath{\mat{K}}}

\providecommand{\mV}{\ensuremath{\mat{V}}}

\providecommand{\vv}{\ensuremath{\vec{v}}}
\providecommand{\vw}{\ensuremath{\vec{w}}}
\providecommand{\vx}{\ensuremath{\vec{x}}}

\newcommand{\bw}{\mathbf{w}}
\newcommand{\bc}{\mathbf{c}}

\numberwithin{equation}{section}
\newtheorem{theorem}{Theorem}[section]

\newtheorem{proposition}[theorem]{Proposition}

\newtheorem{remark}[theorem]{Remark}

\begin{document}
	\title{An asymptotic-preserving dynamical low-rank method for the multi-scale multi-dimensional linear transport equation\footnote{JH's research was supported in part by NSF grant DMS-1620250 and NSF CAREER grant DMS-1654152.}}
	
    \author{Lukas Einkemmer\footnote{Department of Mathematics, University of Innsbruck, Innsbruck, A-6020, Austria (lukas.einkemmer@uibk.ac.at).}, \
	    Jingwei Hu\footnote{Department of Mathematics, Purdue University, West Lafayette, IN 47907, USA (jingweihu@purdue.edu).}, \  
	    and \ Yubo Wang\footnote{Department of Mathematics, Purdue University, West Lafayette, IN 47907, USA (wang3158@purdue.edu).}}   
	\maketitle
	
\begin{abstract}
    We propose a dynamical low-rank method to reduce the computational complexity for solving the multi-scale multi-dimensional linear transport equation. The method is based on a macro-micro decomposition of the equation; the low-rank approximation is only used for the micro part of the solution. The time and spatial discretizations are done properly so that the overall scheme is second-order accurate (in both the fully kinetic and the limit regime) and asymptotic-preserving (AP). That is, in the diffusive regime, the scheme becomes a macroscopic solver for the limiting diffusion equation that automatically captures the low-rank structure of the solution. Moreover, the method can be implemented in a fully explicit way and is thus significantly more efficient compared to the previous state of the art. We demonstrate the accuracy and efficiency of the proposed low-rank method by a number of four-dimensional (two dimensions in physical space and two dimensions in velocity space) simulations.
\end{abstract}

{\small 

}


\section{Introduction}

The linear transport equation models particles such as neutrons or photons interacting with a background medium. This integro-differential equation is widely used in many science and engineering disciplines \cite{Chandrasekhar, Davison}. The linear transport equation belongs to the class of kinetic equations and is consequently posed in a five-dimensional phase space (3D in physical variable and 2D in angle or normalized velocity). This, in particular, implies that its full numerical simulation can be extremely expensive. The situation is further complicated if the scattering strength varies in space by several orders of magnitude, i.e., the equation is stiff in certain region of the domain and non-stiff elsewhere, then an explicit numerical scheme must resolve the smallest collision length.

Recently, a class of dynamical low rank methods has been introduced to solve high dimensional kinetic equations such as the Vlasov equation \cite{El18, EL18_cons, piazzola}, the Boltzmann--BGK equation \cite{E18}, and radiation transport equations \cite{peng2019low,DEL19}. Motivated by these advances, we develop an efficient numerical method to solve the multi-scale multi-dimensional linear transport equation. Employing a low-rank approximation is particularly relevant for the linear transport equation as in the kinetic regime a fine spatial and angular resolution is often required in practice. The low-rank integrator then reduces the five-dimensional problem to a set of (at most) three dimensional equations and thus results in a drastic reduction of memory as well as increased computational efficiency.

An additional goal in the present work is to capture the corresponding asymptotic limit. Although it has been shown in \cite{DEL19} that this can, in principle, be achieved within a low-rank approximation, it comes at the cost of a fully implicit scheme.
The key difference from previous works lies therefore in that, instead of applying the low rank approximation to the unknown distribution function directly, we start with a macro-micro decomposition of the equation and apply the low rank method only to the micro part of the solution. This approach naturally captures the diffusion limit using a more efficient implicit-explicit (IMEX) discretization strategy. In addition, the micro part of the solution becomes low rank in the diffusion limit, hence the method is particularly efficient in this regime.

We mention that the design of numerical schemes that are consistent with certain asymptotic limits falls into the general umbrella of the so-called asymptotic-preserving (AP) schemes \cite{Jin99}, which have been developed for various kinds of kinetic and hyperbolic equations in the past decades, see \cite{Jin_Rev, DD17, HJL17} for an overview. In particular, for the linear transport equation, the use of macro-micro decomposition to achieve the AP property in the diffusive regime first appeared in \cite{LM08}. The stability of the scheme was proved in \cite{LM10} using energy estimates. Comparing to \cite{LM08}, the new difficulty arising in the context of the dynamical low-rank method is to justify the asymptotic limit under the additional projection operator splitting, which we carefully study in this paper. Furthermore, the usual way to generalize the first order (in time) scheme to high order using IMEX Runge-Kutta (RK) schemes, as in \cite{BPR13, JLQX15}, cannot be applied to the low rank case again due to the operator splitting. Hence another contribution of this work is to propose an AP dynamical low-rank method that remains second order in both kinetic and diffusive regimes.

The rest of this paper is organized as follows. In Section~\ref{sec:lineartransport}, we briefly describe the linear transport equation and its macro-micro decomposition. Section~\ref{sec:lowrank} is the main part of the paper where we introduce the dynamical low-rank method. Both the first and second order schemes along with their AP property are discussed in detail. Section~\ref{sec:Fourier} provides a simple Fourier analysis for the solution to the linear transport equation. Section~\ref{sec:num} presents several numerical tests for the two-dimensional equation, where we carefully examine the accuracy, efficiency, rank dependence, and AP property of the proposed method. The paper is concluded in Section~\ref{sec:con}.

\section{The linear transport equation and its macro-micro decomposition}
\label{sec:lineartransport}

We are interested in the following time dependent linear transport equation in diffusive scaling:
\be
\p_tf+\fl 1 \varepsilon \vv\cdot \nabla_{\vx} f = \fl {\sigma^S}{\varepsilon^2}\left(\fl 1 {4\pi}\langle  f\rangle_{\vv}-f\right)-\sigma^Af+G,
\label{eqn:3dlineartransport}
\ee
where $f=f(t,\vx,\vv)$ is the distribution function of time $t$, position $\vx = (x,y,z) \in \Omega_{\vx}\subset\mathbb{R}^3$, and velocity $\vv = (\xi,\eta,\gamma) \in \mathbb{S}^2$ which is confined to the unit sphere\footnote{In the context of radiative transfer, $\vv$ is usually referred to as angle or direction.}. $\langle \ \rangle_{\vv}$ denotes the integration over $\mathbb{S}^2$ with respect to $\vv$. $\sigma^S(\vx)\geq \sigma^S_{\text{min}}>0$ and $\sigma^A(\vx)\geq 0$ are the scattering and absorption coefficients, and $G(\vx)$ is a given source term. Finally $\varepsilon$ is the rescaled collision length, which can range between the kinetic regime $\varepsilon\sim O(1)$ and the diffusive regime $\varepsilon\ll 1$.

The density $\rho = \fl 1 {4\pi}\langle f\rangle_{\vv}$ is defined as the angular average of $f$. In the limit $\varepsilon\rightarrow 0$, $\rho$ satisfies a diffusion equation which can be seen via the Chapman-Enskog expansion. Indeed, (\ref{eqn:3dlineartransport}) can be written as
\be
f = \rho-\varepsilon\fl 1 {\sigma^S} \vv\cdot \nabla_{\vx} f-\varepsilon^2\fl 1 {\sigma^S}\left(\p_tf+\sigma^Af-G\right) = \rho-\varepsilon\fl 1 {\sigma^S} \vv\cdot \nabla_{\vx} \rho + O(\varepsilon^2).
\label{eqn:chap}
\ee
On the other hand, taking $\frac{1}{4\pi}\langle \ \rangle_{\vv}$ of (\ref{eqn:3dlineartransport}) yields
\be
\p_t\rho+\fl 1 {4\pi\varepsilon}\nabla_{\vx} \cdot \langle \vv f\rangle_{\vv}=-\sigma^A\rho+G,
\ee
which, upon substitution of (\ref{eqn:chap}), becomes
\be
\p_t\rho-\nabla_{\vx}\cdot\left(D\nabla_{\vx}\rho\right)=-\sigma^A\rho+G+O(\varepsilon),
\ee
with the diffusion matrix $D$ given by
\be
D = \fl {1} {4\pi\sigma^S} \langle \vv\otimes \vv\rangle_{\vv}= \fl 1 {3\sigma^S} I_{3 \times 3}.
\ee
Therefore, as $\varepsilon \rightarrow 0$ the limit of (\ref{eqn:3dlineartransport}) is the diffusion equation
\be
\p_t\rho-\nabla_{\vx}\cdot\left(\frac{1}{3\sigma^S}\nabla_{\vx}\rho\right)=-\sigma^A\rho+G.
\label{eqn:limiting1}
\ee

In the macro-micro decomposition \cite{LM08}, we write $f$ as
\be \label{ff}
f(t,\vx,\vv) = \rho(t,\vx)+\varepsilon g(t,\vx,\vv),
\ee
where $\rho$ is the macro part of the solution and $g$ is the micro part. Note that $\langle g\rangle_{\vv} = 0$. Substituting (\ref{ff}) into (\ref{eqn:3dlineartransport}) and taking $\frac{1}{4\pi}\langle \ \rangle_{\vv}$, we obtain
\be
\bal
&\p_t\rho +\fl 1 {4\pi} \nabla_{\vx} \cdot \langle \vv  g\rangle_{\vv}=-\sigma^A\rho+G.
\eal
\label{eqn:reducedmm1}
\ee
Subtracting (\ref{eqn:reducedmm1}) from (\ref{eqn:3dlineartransport}) yields
\be
\p_tg+\fl 1 \varepsilon\left(I-\fl 1 {4\pi}\langle \ \rangle_{\vv}\right)\left(\vv \cdot \nabla_{\vx} g\right)+\fl 1 {\varepsilon^2} \vv \cdot \nabla_{\vx} \rho=-\fl {\sigma^S} {\varepsilon^2} g-\sigma^Ag.
\label{eqn:reducedmm2}
\ee

The coupled system (\ref{eqn:reducedmm1}) and (\ref{eqn:reducedmm2}) is the macro-micro decomposition of the linear transport equation (\ref{eqn:3dlineartransport}). In the limit $\varepsilon \rightarrow 0$, we have from (\ref{eqn:reducedmm2}):
\be
	g = -\fl 1 {\sigma^S}\vv \cdot \nabla_{\vx} \rho,
\ee
which, when substituting into (\ref{eqn:reducedmm1}), yields the same diffusion equation (\ref{eqn:limiting1}).

\section{The dynamical low-rank method for the linear transport equation}
\label{sec:lowrank}

We first constrain $g(t,\vx,\vv)$ to a low rank manifold $\mathcal{M}$ such that
\be \label{lowrankg}
	g(t,\vx,\vv) = \sum_{i,j=1}^rX_i(t,\vx)S_{ij}(t)V_j(t,\vv), 
\ee
where $r$ is called the rank and the basis functions $\{X_i\}_{1\leq i\leq r}$ and $\{V_j\}_{1\leq j \leq r}$ are orthonormal:
\be
{\langle X_i,X_k\rangle}_{\vx}=\delta_{ik}, \ {\langle V_j,V_k\rangle}_{\vv}=\delta_{jk},
\ee
with $\langle \cdot,\cdot\rangle_{\vx}$ and $\langle \cdot,\cdot\rangle_{\vv}$ being the inner products on $L^2(\Omega_{\vx})$ and $L^2(\mathbb{S}^2)$, respectively. 

With this low rank approximation, (\ref{eqn:reducedmm1}) becomes
\be
\label{steprho}
\p_t\rho =-\fl 1 {4\pi}\sum\limits_{i,j=1}^r\nabla_{\vx} \cdot  \left( X_iS_{ij} \langle\vv V_j\rangle_{\vv}  \right)-\sigma^A\rho+G.
\ee
For (\ref{eqn:reducedmm2}), we write
\be \label{eqn:RHS}
\p_tg=-\fl 1 \varepsilon\left(I-\fl 1 {4\pi}\langle \ \rangle_{\vv}\right)\left(\vv \cdot \nabla_{\vx} g\right)-\fl 1 {\varepsilon^2}\vv \cdot \nabla_{\vx} \rho-\fl {\sigma^S} {\varepsilon^2} g-\sigma^Ag :=\text{RHS}.
\ee
Equation \eqref{eqn:RHS}, however, does not uniquely specify the dynamics of the low-rank factors $X_i$, $S_{ij}$, and $V_j$. We therefore impose the following gauge conditions \cite{KL07}:
\be \label{eq:gauge}
{\langle \partial_t X_i,X_k\rangle}_{\vx}=0, \ {\langle \partial_t V_j,V_k\rangle}_{\vv}=0.
\ee
Let us emphasize that the resulting dynamics of $g$ is independent of the specific gauge conditions chosen. However, using \eqref{eq:gauge} is convenient as it allows us to easily obtain evolution equations in terms of the low-rank factors. To that end, we now project the right hand side of (\ref{eqn:RHS}) onto the tangent space of $\mathcal{M}$:
\be
	\p_t g = P_g(\text{RHS}),
\label{projected equation}
\ee
where the orthogonal projector $P_g$ can be written as
\be
	P_g(\text{RHS}) = \sum_{j=1}^r\langle V_j,\text{RHS}\rangle_{\vv} V_j - \sum_{i,j=1}^rX_i\langle X_iV_j,\text{RHS}\rangle_{\vx,\vv}V_j+\sum_{i=1}^rX_i\langle X_i,\text{RHS}\rangle_{\vx}. \label{projected_parts}
\ee
Using \eqref{projected_parts} and the gauge conditions we can in principle derive evolution equations for $X_i$, $S_{ij}$, and $V_j$. However, this process requires inverting the matrix $S=(S_{ij})$. Since an accurate approximation mandates that $S$ has small singular values, the resulting problem is severely ill-conditioned. Thus, we will use the projector splitting scheme introduced in \cite{LO14}. For a corresponding mathematical analysis see \cite{KLW16}. This scheme has been extensively used in the literature, see e.g.~\cite{El18,peng2019low,Lubich15tii}, and extensions to various tensor formats have also been proposed \cite{LVW2018,LOV2015}. The main idea is to split equation \eqref{projected equation} into the following three subflows
\begin{align*}
    \partial_t g &= \sum_{j=1}^r\langle V_j,\text{RHS}\rangle_{\vv} V_j, \\
    \partial_t g &= -\sum_{i,j=1}^rX_i\langle X_iV_j,\text{RHS}\rangle_{\vx,\vv}V_j, \\
    \partial_t g &= \sum_{i=1}^rX_i\langle X_i,\text{RHS}\rangle_{\vx}.
\end{align*}
This is particularly convenient as for the first subflow $V_j$ is constant (in time), for the third subflow $X_i$ is constant, and for the second subflow both $X_i$ and $V_j$ are constant. Thus, we can write
\begin{align}
    \partial_t K_j  &= \langle V_j,\text{RHS}\rangle_{\vv}, \label{eq:K-step}\\
    \partial_t S_{ij} &= -\langle X_iV_j,\text{RHS}\rangle_{\vx,\vv}, \label{eq:S-step}\\
    \partial_t L_i  &= \langle X_i,\text{RHS}\rangle_{\vx}, \label{eq:L-step}
\end{align}
where
\be
K_j(t,\vx)=\sum\limits_{i=1}^rX^{}_{i}(t,\vx)S^{}_{ij}(t), \qquad\qquad L_i(t,\vv)=\sum\limits_{j=1}^rS^{}_{ij}(t)V^{}_{j}(t,\vv).
\ee
After solving each subflow we use a QR decomposition to obtain $X_i$ and $S_{ij}$ from $K_j$ and $S_{ij}$ and $V_j$ from $L_i$, respectively.

\subsection{A first order in time scheme}
\label{subsec:first}

Our goal is to solve the coupled system (\ref{steprho}) and (\ref{projected equation}) using the projector splitting integrator outlined in the previous section. We now proceed by deriving the evolution equations corresponding to the subflows given by equations \eqref{eq:K-step}-\eqref{eq:L-step}.
\begin{itemize}
    \item {\bf $K$-step}: Solve $\p_t K_j = \langle V_j,\text{RHS}\rangle_{\vv}$ with $\{V_j\}_{1\leq j \leq r}$ unchanged.
\be
\bal
\label{stepK}
\p_t K_j =& \langle V_j,\text{RHS}\rangle_{\vv}\\
=&-\fl 1 \varepsilon \sum_{l=1}^r \left(\langle \vv V_j^{}V_l^{}\rangle_{\vv}-\fl 1 {4\pi}\langle V_j^{}\rangle_{\vv}\langle \vv V_l^{}\rangle_{\vv}\right)\cdot \nabla_{\vx} K_l \\
&-\fl 1 {\varepsilon^2} \langle\vv V_j\rangle_{\vv}\cdot\nabla_{\vx}\rho-\left(\fl {\sigma^S}{\varepsilon^2}+\sigma^A\right)K_j.
\eal
\ee
\item {\bf $L$-step}: Solve $\p_tL_i = \langle X_i,\text{RHS}\rangle_{\vx}$ with $\{X_i\}_{1\leq i\leq r}$ unchanged. 
\be
\bal
\label{stepL}
\p_tL_i =& \langle X_i,\text{RHS}\rangle_{\vx}\\
=&-\fl 1 \varepsilon \sum_{k=1}^r  \left( \vv L_k-\frac{1}{4\pi}\langle\vv L_k\rangle_{\vv}\right)\cdot \langle X_i^{}\nabla_{\vx}X_k^{}\rangle_{\vx} \\
&-\fl 1 {\varepsilon^2} \vv\cdot\langle X_i \nabla_{\vx}\rho\rangle_{\vx}-\sum_{k=1}^r\left\langle X_i^{}\left(\fl {\sigma^S}{\varepsilon^2}+\sigma^A\right)X_k^{}\right\rangle_{\vx} L_k.
\eal
\ee
\item {\bf $S$-step}: Solve $\p_tS_{ij} =-\langle X_iV_j,\text{RHS}\rangle_{\vx,\vv}$ with both $\{X_i\}_{1\leq i\leq r}$ and $\{V_j\}_{1\leq j \leq r}$ unchanged. 
\be
\label{stepS}
\bal
 \p_tS_{ij} =&-\langle X_iV_j,\text{RHS}\rangle_{\vx,\vv}\\
=& \fl 1 {\varepsilon} \sum_{k,l=1}^r \left(\langle \vv V_j^{}V_l^{}\rangle_{\vv}-\fl 1 {4\pi}\langle V_j^{}\rangle_{\vv}\langle \vv V_l^{}\rangle_{\vv}\right)\cdot \langle X_i\nabla_{\vx}X_k\rangle_{\vx} S_{kl}\\
&+\fl 1 {\varepsilon^2} \langle\vv V_j\rangle_{\vv}\cdot \langle X_i\nabla_{\vx}\rho\rangle_{\vx}+\sum_{k=1}^r\left\langle X_i^{}\left(\fl {\sigma^S}{\varepsilon^2}+\sigma^A\right)X_k^{}\right\rangle_{\vx} S_{kj}.
\eal
\ee
\end{itemize}

Therefore, for the overall system, we can construct a simple first order in time scheme. Suppose at time step $t^n$, we have $(X_i^n,V_j^n,S_{ij}^n,\rho^n)$. To obtain the solution $(X_i^{n+1},V_j^{n+1},S_{ij}^{n+1},\rho^{n+1})$ at $t^{n+1}$ we proceed as follows:
\begin{enumerate}
\item {\bf $K$-step}: Solve (\ref{stepK}) for a full time step $\dt t$, update from $(X_i^n,V_j^n,S_{ij}^n)$ to $(X_i^{n+1},V_j^n,S_{ij}^{(1)})$ using $\rho^n$. Specifically, given $K_j^n= \sum\limits_{i=1}^rX^{n}_{i}S^{n}_{ij}$, we discretize (\ref{stepK}) using a first order IMEX scheme (i.e., forward-backward Euler scheme) as
\be
\label{semistepK}
\bal
	\fl {K_j^{n+1}-K_j^{n}}{\dt t} =& -\fl 1 \varepsilon \sum_{l=1}^r \left(\langle \vv V_j^{n}V_l^{n}\rangle_{\vv}-\fl 1 {4\pi}\langle V_j^{n}\rangle_{\vv}\langle \vv V_l^{n}\rangle_{\vv}\right)\cdot \nabla_{\vx} K_l^n \\
        &-\fl 1 {\varepsilon^2} \left( \langle\vv V_j^n\rangle_{\vv}\cdot\nabla_{\vx}\rho^{n}+ {\sigma^S}K_j^{n+1}\right)-\sigma^AK_j^{n},
\eal
\ee
where the term $\sigma^SK_j$ is treated implicitly to overcome the stiffness induced by a small $\varepsilon$. We then perform the QR decomposition of $K_j^{n+1}$ to obtain the updated basis functions $X_i^{n+1}$ and the matrix $S_{ij}^{(1)}$:
\be
K_j^{n+1}=\sum_{i=1}^rX_i^{n+1}S_{ij}^{(1)}.
\ee

\item {\bf $L$-step}: Solve (\ref{stepL}) for a full time step $\dt t$, update from $(X_i^{n+1},V_j^n,S_{ij}^{(1)})$ to $(X_i^{n+1},V_j^{n+1},S_{ij}^{(2)})$ using $\rho^n$. Specifically, given $L_i^n=\sum\limits_{j=1}^rS_{ij}^{(1)}V_{j}^n$, we discretize (\ref{stepL}) (similar to (\ref{stepK})) as follows
\be
\label{semistepL}
\bal
\fl {L_i^{n+1}-L_i^n}{\dt t} =& -\fl 1 \varepsilon \sum_{k=1}^r  \left( \vv L_k^n-\frac{1}{4\pi}\langle\vv L_k^n\rangle_{\vv}\right) \cdot \langle X_i^{n+1}\nabla_{\vx}X_k^{n+1}\rangle_{\vx} \\
        & -\fl 1 {\varepsilon^2} \left( \vv\cdot\langle X_i^{n+1} \nabla_{\vx}\rho^{n}\rangle_{\vx} +\sum_{k=1}^r\left\langle X_i^{n+1}\sigma^S X_k^{n+1}\right\rangle_{\vx} L_k^{n+1}\right) -\sum_{k=1}^r\left\langle X_i^{n+1}\sigma^A X_k^{n+1}\right\rangle_{\vx} L_k^{n}.
\eal
\ee
We then perform the QR decomposition of $L_i^{n+1}$ to obtain the updated basis $V_{j}^{n+1}$ and matrix $S_{ij}^{(2)}$:
\be
L_i^{n+1}=\sum\limits_{j=1}^rS_{ij}^{(2)}V_{j}^{n+1}.
\ee

\item {\bf $S$-step}: Solve (\ref{stepS}) for a full time step $\dt t$, update from $(X_i^{n+1},V_j^{n+1},S_{ij}^{(2)})$ to $(X_i^{n+1},V_j^{n+1},S_{ij}^{n+1})$ using $\rho^n$. Specifically, given $S_{ij}^{(2)}$, we discretize (\ref{stepS}) (similar to (\ref{stepK})) as follows
\be
\label{semistepS}
\bal
\fl {S^{n+1}_{ij}-S^{(2)}_{ij}}{\dt t} =& \fl 1 {\varepsilon} \sum_{k,l=1}^r\left(\langle \vv V^{n+1}_jV^{n+1}_l\rangle_{\vv}-\fl 1 {4\pi}\langle V_j^{n+1}\rangle_{\vv}\langle \vv V_l^{n+1}\rangle_{\vv}\right)\cdot\langle X^{n+1}_i\nabla_{\vx}X^{n+1}_k\rangle_{\vx}S_{kl}^{(2)}
\\&+\fl 1 {\varepsilon^2}\left(\langle\vv V^{n+1}_j\rangle_{\vv}\cdot  \langle X_i^{n+1}\nabla_{\vx}\rho^{n}\rangle_{\vx}+\sum_{k=1}^r \langle X_i^{n+1}\sigma^SX_k^{n+1}\rangle_{\vx} S_{kj}^{n+1}\right)\\
&+\sum_{k=1}^r \langle X_i^{n+1}\sigma^AX_k^{n+1}\rangle_{\vx} S_{kj}^{(2)}.
\eal
\ee

\item {\bf $\rho$-step}: Solve (\ref{steprho}) for a full time step $\dt t$, update from $\rho^{n}$ to $\rho^{n+1}$ using $(X_i^{n+1},V_j^{n+1},S_{ij}^{n+1})$. Specifically, given $\rho^n$, we discretize (\ref{steprho}) as
\be
\label{semisteprho}
\fl {\rho^{n+1}-\rho^{n}}{\dt t} =-\fl 1 {4\pi}\sum\limits_{i,j=1}^r \nabla_{\vx} \cdot  \left( X_i^{n+1}S_{ij}^{n+1} \langle\vv V_j^{n+1}\rangle_{\vv}  \right) -\sigma^A\rho^{n}+G.
\ee
\end{enumerate}

For clarity, we will refer to the above scheme as the {\bf $K$-$L$-$S$-$\rho$ scheme} in the following.

\subsection{AP property of the first order scheme}
\label{sec:firstAP}

In this subsection, we analyze the AP property of the first order scheme introduced in the previous section. Our conclusion is summarized in the following proposition.

\begin{proposition} \label{prop1}
    In the limit $\varepsilon\rightarrow 0$, the first order IMEX $K$-$L$-$S$-$\rho$ scheme (i.e., (\ref{semistepK}), (\ref{semistepL}), (\ref{semistepS}), and (\ref{semisteprho})) becomes the forward Euler scheme for the limiting diffusion equation (\ref{eqn:limiting1}), provided  that for the initial value we have $\xi, \eta, \gamma \in \text{span}(\{V_j^0\}_{j=1}^r)$. 
\end{proposition}
\begin{remark}
    If, for a given initial value $(X_i^0, S_{ij}^0, V_j^0)$, one of the conditions $\xi, \eta, \gamma \in \text{span}(\{V_j^0\}_{j=1}^r)$ is not satisfied, we can simply add them to the approximation space. For example, if $\xi \not\in \text{span}(\{V_j^0\}_{j=1}^r)$, we consider 
    \[ \tilde{X}^0  = [X_1^0, \dots, X_r^0, h], \qquad \tilde{S}^0 = \bmat{S^0 & 0\\ 0 & 0}, \qquad \tilde{V}^0 = [V_1^0, \dots, V_r^0, \xi], \]
    where $h$ is an arbitrary function. We then orthogonalize $\tilde{X}^0$ and $\tilde{V}^0$ (e.g.~using the Gram-Schmidt process) and use the result as the initial value in our algorithm. This increases the rank to at most $r+3$.

\end{remark}

\begin{proof}
    In the {\bf K-step}, let $\varepsilon \rightarrow 0$, we have from (\ref{semistepK}): 
\be \label{KK}
K_j^{n+1}=-\langle \vv V_j^n\rangle_{\vv}\cdot\fl {\nabla_{\vx}\rho^n}{\sigma^S}.
\ee
Without loss of generality, we assume that the three components of $\fl {\nabla_{\vx}\rho^n}{\sigma^S}$: $\fl {\p_x \rho^{n}}{\sigma^S}$, $\fl {\p_y \rho^{n}}{\sigma^S}$ and $\fl {\p_z \rho^{n}}{\sigma^S}$ are linearly independent\footnote{If they are linearly dependent, say, $\text{span}\{\fl {\p_x \rho^{n}}{\sigma^S}, \fl {\p_y \rho^{n}}{\sigma^S}, \fl {\p_z \rho^{n}}{\sigma^S}\}=\text{span}\{\fl {\p_x \rho^{n}}{\sigma^S}\}$, then one just needs to replace the second and third components of $X_0$ by $X_2^{n+1}$ and $X_3^{n+1}$ and the same analysis carries over.}. Then after the QR decomposition of $K_j^{n+1}$, the span of the new basis $\{X_i^{n+1}\}_{1\leq i\leq 3}$ would be the same as $\text{span}\{\fl {\p_x \rho^{n}}{\sigma^S}, \fl {\p_y \rho^{n}}{\sigma^S}, \fl {\p_z \rho^{n}}{\sigma^S}\}$. In other words, we can write
\be \label{XX}
X^{n+1} :=\bmat{ X_1^{n+1} & X_2^{n+1} & X_3^{n+1} & X_4^{n+1} & \cdots &X_r^{n+1} }= \underbrace{\bmat{\fl {\p_x \rho^{n}}{\sigma^S}&\fl {\p_y \rho^{n}}{\sigma^S} &\fl {\p_z \rho^n}{\sigma^S}&  X_4^{n+1} & \cdots &X_r^{n+1}}}_{:=X_0}D_1,
\ee
where $D_1$ is an invertible $r \times r$ matrix. 

In the {\bf L-step}, let $\varepsilon \rightarrow 0$, we have from (\ref{semistepL}):
\be
	\sum_{k=1}^r\langle X_i^{n+1}\sigma^S X^{n+1}_k\rangle_{\vx} L_k^{n+1} = -\vv \cdot \langle X_i^{n+1}\nabla_{\vx}\rho^n\rangle_{\vx}.
\ee
    Since the matrix $A:=(\langle X_i^{n+1}\sigma^S X^{n+1}_k\rangle_{\vx})_{1\leq i\leq r, 1\leq k \leq r}$ is symmetric positive definite (since $\sigma^S>0$), hence invertible (whose inverse, say, is matrix $B=(b_{ki})_{1\leq k\leq r, 1\leq i \leq r}$), we have
\be
    L_k^{n+1} = - \vv \cdot \left( \sum_{i=1}^r b_{ki}\langle X_i^{n+1}\nabla_{\vx}\rho^n\rangle_{\vx} \right).
\ee
After the QR decomposition of $L_k^{n+1}$, we can write (by a similar argument as above)
\be \label{VV}
V^{n+1}:=\bmat{ V_1^{n+1} & V_2^{n+1} & V_3^{n+1} & V_4^{n+1} & \cdots &V_r^{n+1} }= \underbrace{\bmat{\xi & \eta & \gamma & V_4^{n+1} & \cdots &V_r^{n+1}}}_{:=V_0}D_2,
\ee
where $D_2$ is an invertible $r \times r$ matrix.

In the {\bf S-step}, let $\varepsilon \rightarrow 0$, we have from (\ref{semistepS}):
\be
\bal
	\sum_{k=1}^r\langle X^{n+1}_i\sigma^S X^{n+1}_k\rangle_{\vx} S^{n+1}_{kj}&=-\langle \vv V_j^{n+1}\rangle_{\vv} \cdot \langle X_i^{n+1}\nabla_{\vx}\rho^n\rangle_{\vx} \\
	&= -\langle X^{n+1}_i\vv \cdot \nabla_{\vx}\rho^nV_j^{n+1}\rangle_{\vx,\vv}.
\eal
\label{eqn:updateS}
\ee
We may write (\ref{eqn:updateS}) as $AS^{n+1}=C$. Since the matrix $A$ is invertible, we know that the matrix $S^{n+1}$ is unique.
We next claim that the $S^{n+1}$ defined as
\be
S^{n+1} := D_1^{-1}\bmat{-I_{3\times 3}&0\\ 0 & 0}D_2^{-T},
\ee
satisfies (\ref{eqn:updateS}), where the middle matrix is of size $r\times r$, with $-I_{3\times 3}$ in the first $3\times 3$ block and zero elsewhere. Indeed, using (\ref{XX}) and (\ref{VV}) we have
\be
\bal
	g^{n+1}=\sum_{i,j=1}^rX_i^{n+1}S_{ij}^{n+1}V_j^{n+1}=X^{n+1}S^{n+1}(V^{n+1})^T = X_0\bmat{-I_{3\times 3}&0\\ 0 & 0}V_0^T=-\vv \cdot  \frac{\nabla_{\vx}\rho^n}{\sigma^S}.
\eal 
\label{eqn:limiteqn}
\ee
Therefore,
\be
	(X^{n+1})^T\sigma^SX^{n+1}S^{n+1}(V^{n+1})^TV^{n+1} = -(X^{n+1})^T(\vv \cdot  \nabla_{\vx}\rho^n)V^{n+1},
\ee
which, upon taking $\langle \ \rangle_{\vx,\vv}$, yields
\be
\langle (X^{n+1})^T\sigma^S X^{n+1}\rangle_{\vx}  S^{n+1} = -\langle (X^{n+1})^T(\vv \cdot  \nabla_{\vx}\rho^n)V^{n+1} \rangle_{\vx,\vv},
\ee
which is precisely (\ref{eqn:updateS}). 

On the other hand, substituting (\ref{eqn:limiteqn}) into (\ref{semisteprho}) gives
\be \label{forwardEuler}
\fl {\rho^{n+1}-\rho^{n}}{\dt t} =\nabla_{\vx} \cdot  \left(\frac{1}{3\sigma^S} \nabla_{\vx}\rho^n\right)-\sigma^A\rho^n+G,
\ee
which is the forward Euler scheme for the limiting diffusion equation (\ref{eqn:limiting1}).
\end{proof}

\subsection{Some other first order schemes and their AP property}
\label{sec:otherAP}

From the operator splitting point of view, the previously introduced $K$-$L$-$S$-$\rho$ scheme is certainly not the only first order scheme. In fact, one can switch the order of $K$, $L$, and $S$ steps arbitrarily and still obtains a first order scheme. For example, the $L$-$K$-$S$-$\rho$ scheme is also first order and preserves the same asymptotic limit as the $K$-$L$-$S$-$\rho$ scheme (since the proof of Proposition~\ref{prop1} still holds if one switches the $K$ and $L$ steps). Nonetheless, for some other first order schemes, such as $L$-$S$-$K$-$\rho$, $S$-$L$-$K$-$\rho$, $K$-$S$-$L$-$\rho$, and $S$-$K$-$L$-$\rho$ schemes, their AP property needs to be examined individually. Fortunately, as we will show in the following, by slightly different arguments these schemes all have the same asymptotic limit as the $K$-$L$-$S$-$\rho$ scheme.

\begin{itemize}
\item {\bf $L$-$S$-$K$-$\rho$ scheme} and {\bf $S$-$L$-$K$-$\rho$ scheme}. 

After the first two substeps ($L$-$S$ or $S$-$L$), the span of the updated basis $\{V_j^{n+1}\}_{1\leq j \leq r}$ will contain $\vv$. After the substep $K$, one has $K_j^{n+1}=-\langle \vv V_j^{n+1}\rangle_{\vv}\cdot\fl {\nabla_{\vx}\rho^n}{\sigma^S}$. Hence,
\be
g^{n+1}=\sum\limits_{j=1}^r K_j^{n+1}V_j^{n+1}=-\sum_{j=1}^r \langle \vv V_j^{n+1}\rangle_{\vv} V_j^{n+1}\cdot \fl {\nabla_{\vx}\rho^n}{\sigma^S}=-\vv \cdot \fl {\nabla_{\vx}\rho^n}{\sigma^S}.
\ee
Substituting $g^{n+1}$ into the last $\rho$ step recovers (\ref{forwardEuler}).

\item {\bf $K$-$S$-$L$-$\rho$ scheme} and {\bf $S$-$K$-$L$-$\rho$ scheme}.

After the first two substeps ($K$-$S$ or $S$-$K$), one has 
\be 
\bmat{ X_1^{n+1} &  \cdots &X_r^{n+1} }=\bmat{\fl {\p_x \rho^{n}}{\sigma^S}&\fl {\p_y \rho^{n}}{\sigma^S} &\fl {\p_z \rho^n}{\sigma^S}&  X_4^{n+1} & \cdots &X_r^{n+1}}D_1,
\ee
where $D_1$ is an invertible $r \times r$ matrix. After the substep $L$, one has
\be \label{LL}
	\sum_{k=1}^r\langle X_i^{n+1}\sigma^S X^{n+1}_k\rangle_{\vx} L_k^{n+1} = -\vv \cdot \langle X_i^{n+1}\nabla_{\vx}\rho^n\rangle_{\vx},
\ee
and $\{L_k^{n+1}\}_{1\leq k \leq r}$ is uniquely determined since the matrix $\langle X_i^{n+1}\sigma^S X^{n+1}_k\rangle_{\vx}$ is invertible. We now claim that $\{L_k^{n+1}\}_{1\leq k \leq r}$ defined as follows
\be
\bmat{ L_1^{n+1} &  \cdots &L_r^{n+1} }:=-\bmat{\xi & \eta & \gamma & 0 & \cdots &0}D_1^{-T}.
\ee
satisfies (\ref{LL}). Indeed, for such $L_k$, one has
\be
g^{n+1}=\sum_{k=1}^r X_k^{n+1}L_k^{n+1}=-\vv \cdot \fl {\nabla_{\vx}\rho^n}{\sigma^S} \  \Longrightarrow \sum_{k=1}^r \sigma^S X_k^{n+1}L_k^{n+1} =-\vv \cdot \nabla_{\vx}\rho^n,
\ee
which, upon projection onto the space spanned by $\{X_i^{n+1}\}_{1\leq i \leq r}$, yields (\ref{LL}). On the other hand, substituting $g^{n+1}$ into the last $\rho$ step recovers (\ref{forwardEuler}).
\end{itemize}

\begin{remark}
The discussion in this subsection implies that one has the flexibility to choose the updating order of $K$, $L$ and $S$, while still maintaining the AP property. This flexibility is crucial in designing second order schemes, where one needs to properly compose these steps to achieve high order as well as preserve the asymptotic limit.
\end{remark}

\subsection{A second order in time scheme and its AP property}
\label{subsec:second}

We now extend the first order scheme to second order. Due to the operator splitting necessary in the low rank method, a straightforward application of the IMEX-RK scheme as used in \cite{BPR13, JLQX15} does not work (there a coupled system for $\rho$ and $g$ is solved simultaneously; in the present work $\rho$ has to be ``frozen" while updating $g$). In the following, we propose a scheme that maintains second order in both kinetic and diffusive regimes. It is a proper combination of the almost symmetric Strang splitting \cite{EO14,le2014almost} and the IMEX-RK scheme.

Suppose at time step $t^n$, we have $(X_i^n,V_j^n,S_{ij}^n,\rho^n)$. To obtain the solution $(X_i^{n+1},V_j^{n+1},S_{ij}^{n+1},\rho^{n+1})$ at $t^{n+1}$, we proceed as follows:
\begin{enumerate}
\item {\bf $\rho$-step}: Solve (\ref{steprho}) for a half time step $\dt t/2$, update from $\rho^n$ to $\rho^{n+\fl 1 2}$ using $(X_i^{n},V_j^{n},S_{ij}^{n})$. 
\item {\bf $K$-step}: Solve (\ref{stepK}) for a half time step $\dt t/2$, update from $(X_i^n,V_j^n,S_{ij}^n)$ to $(X_i^{n+\fl 1 2},V_j^n,S_{ij}^{(1)})$ using $\rho^{n+\fl 1 2}$.
\item {\bf $L$-step}: Solve (\ref{stepL}) for a half time step $\dt t/2$, update from $(X_i^{n+\fl 1 2},V_j^n,S_{ij}^{(1)})$ to $(X_i^{n+\fl 1 2},V_j^{n+\fl 1 2},S_{ij}^{(2)})$ using $\rho^{n+\fl 1 2}$.
\item {\bf $S$-step}: Solve (\ref{stepS}) for a half time step $\dt t/2$, update from $(X_i^{n+\fl 1 2},V_j^{n+ \fl 1 2},S_{ij}^{(2)})$ to $(X_i^{n+\fl 1 2},V_j^{n+\fl 1 2},S_{ij}^{n+\fl 1 2})$ using $\rho^{n+ \fl 1 2}$.
\item {\bf $S$-step}: Solve (\ref{stepS}) for a half time step $\dt t/2$, update from $(X_i^{n+\fl 1 2},V_j^{n+ \fl 1 2},S_{ij}^{n+ \fl 1 2})$ to $(X_i^{n+\fl 1 2},V_j^{n+\fl 1 2},S_{ij}^{(3)})$ using $\rho^{n+ \fl 1 2}$.
\item {\bf $L$-step}: Solve (\ref{stepL}) for a half time step $\dt t/2$, update from $(X_i^{n+\fl 1 2},V_j^{n+ \fl 1 2},S_{ij}^{(3)})$ to $(X_i^{n+\fl 1 2},V_j^{n+1},S_{ij}^{(4)})$ using $\rho^{n+\fl 1 2}$.
\item {\bf $K$-step}: Solve (\ref{stepK}) for a half time step $\dt t/2$, update from $(X_i^{n+\fl 1 2},V_j^{n+ 1},S_{ij}^{(4)})$ to $(X_i^{n+1},V_j^{n+1},S_{ij}^{n+1})$ using $\rho^{n+\fl 1 2}$.
\item {\bf $\rho$-step}: Solve (\ref{steprho}) for a full time step $\dt t$, update from $\rho^n$ to $\rho^{n+1}$ using $(X_i^{n+\fl 1 2},V_j^{n+\fl 1 2},S_{ij}^{n+\fl 1 2})$.
\end{enumerate}

More specifically, in step 1, we use the forward Euler scheme to discretize (\ref{steprho}):
\be
\label{semisteprho2_1}
\fl {\rho^{n+\fl 1 2}-\rho^{n}}{\dt t/2} =-\fl 1 {4\pi}\sum\limits_{i,j=1}^r \nabla_{\vx} \cdot  \left( X_i^nS_{ij}^n \langle\vv V_j^n\rangle_{\vv}  \right) -\sigma^A\rho^{n}+G.
\ee
In steps 2-7, we use a second order IMEX-RK scheme to discretize the system for $K$, $L$ or $S$. Let us take step 2 for example,
\be
\bal
	K_j^{(p)}=&K_j^n -\fl {\dt t} 2 \sum_{q= 1}^{p-1}\tilde{a}_{p q}\left(\fl {1} \varepsilon\sum_{l=1}^r \left(\langle \vv V_j^{n}V_l^{n}\rangle_{\vv}-\fl 1 {4\pi}\langle V_j^{n}\rangle_{\vv}\langle \vv V_l^{n}\rangle_{\vv}\right)\cdot \nabla_{\vx} K_l^{(q)}+\fl {1} {\varepsilon^2} \langle\vv V_j^n\rangle_{\vv}\cdot\nabla_{\vx}\rho^{n+\fl 1 2}+\sigma^AK_j^{(q)}\right) \\
	&-\fl {\dt t} 2\sum_{q= 1}^{p}{a}_{pq}\left(\fl {\sigma^S}{\varepsilon^2}K_j^{(q)}\right),\  p= 1,\ldots,s,\\
	K_j^{n+1}=&K_j^n -\fl {\dt t} 2 \sum_{p = 1}^{s}\tilde{w}_{p}\left(\fl {1} \varepsilon\sum_{l=1}^r \left(\langle \vv V_j^{n}V_l^{n}\rangle_{\vv}-\fl 1 {4\pi}\langle V_j^{n}\rangle_{\vv}\langle \vv V_l^{n}\rangle_{\vv}\right)\cdot \nabla_{\vx} K_l^{(p)}+\fl {1} {\varepsilon^2} \langle\vv V_j^n\rangle_{\vv}\cdot\nabla_{\vx}\rho^{n+\fl 1 2}+\sigma^AK_j^{(p)}\right) \\
	&-\fl {\dt t} 2\sum_{p= 1}^{s}{w}_{p}\left(\fl {\sigma^S}{\varepsilon^2}K_j^{(p)}\right),\
\eal
\ee
where $\tilde{A} = (\tilde{a}_{p q})$, $\tilde{a}_{p q} = 0$ for $q\geq p $ and $A = (a_{pq})$, $a_{p q} = 0$ for $q>p$ are $s \times s$ matrices. Along with $\tilde{\vw} = (\tilde{w}_1,\ldots,\tilde{w}_s)^T$, ${\vw} = ({w}_1,\ldots,{w}_s)^T$, they can be represented by a double Butcher tableau:
\begin{equation}
\centering
\begin{tabular}{c|c}
$\tilde{\bc}$ & $\tilde{A}$\\
\hline
& $\tilde{\bw}^T$
\end{tabular} \quad \quad
\begin{tabular}{c|c}
$\bc$ & $A$\\
\hline
& $\bw^T$
\end{tabular}
\end{equation}
where $\tilde{\bc}=(\tilde{c}_1,\dots,\tilde{c}_s)^T$, $\bc=(c_1,\dots,c_s)^T$ are defined as
\begin{equation}
\tilde{c}_{p}=\sum_{q=1}^{p-1}\tilde{a}_{p q}, \quad c_{p}=\sum_{q=1}^{p}a_{p q}.
\end{equation}
Here we employ the ARS(2,2,2) scheme whose double tableau is given by
\begin{equation} 
\centering
\begin{tabular}{c | c c c}
 0& 0 & 0 & 0 \\
 $\gamma$& $\gamma$ & 0 & 0\\
 1& $\delta$ & $1-\delta$ & 0\\
\hline
& $\delta$ &  $1-\delta$ & 0
\end{tabular} \quad \quad
\begin{tabular}{c | c c c}
 0& $0$ & 0 & 0 \\
 $\gamma$& $0$ & $\gamma$ & 0\\
 1& $0$ & $1-\gamma$ & $\gamma$\\
\hline
& $0$ &  $1-\gamma$ & $\gamma$
\end{tabular}
\qquad  \gamma = 1-\fl {\sqrt{2}}{2}, \quad \delta = 1-\fl 1 {2\gamma}.
\end{equation}
Finally, in step 8, we use the midpoint scheme to discretize (\ref{steprho}):
\be
\label{semisteprho2_2}
\fl {\rho^{n+1}-\rho^{n}}{\dt t} =-\fl 1 {4\pi}\sum\limits_{i,j=1}^r \nabla_{\vx} \cdot  \left( X_i^{n+\fl 1 2} S_{ij}^{n + \fl 1 2}\langle\vv V_j^{n+ \fl 1 2}\rangle_{\vv}  \right)  -\sigma^A\rho^{n+\fl 1 2}+G.
\ee

Let us analyze the AP property of the above second order scheme. 
First, steps 2-4 ($K$-$L$-$S$) are (almost) the same as steps 1-3 in the first order $K$-$L$-$S$-$\rho$ scheme (as discussed in Section~\ref{sec:firstAP}), hence as $\varepsilon \rightarrow 0$, one has
\be \label{secondgg1}
g^{n+\fl 12}=\sum_{i,j=1}^rX_i^{n+\fl 1 2}S_{ij}^{n+\fl 1 2}V_j^{n+\fl 1 2} = -\vv \cdot  \frac{\nabla_{\vx}\rho^{n+\fl 1 2}}{\sigma^S}.
\ee

Furthermore, steps 5-6 ($S$-$L$-$K$) are (almost) the same as steps 1-3 in the first order $S$-$L$-$K$-$\rho$ scheme (as discussed in Section~\ref{sec:otherAP}), hence as $\varepsilon \rightarrow 0$, one has
\be \label{secondgg2}
g^{n+1}=\sum_{j=1}^rK_j^{n+1}V_j^{n+1}=-\vv \cdot  \frac{\nabla_{\vx}\rho^{n+\fl 1 2}}{\sigma^S}.
\ee

Finally, substituting (\ref{secondgg2}) into (\ref{semisteprho2_1}) and (\ref{secondgg1}) into (\ref{semisteprho2_2}), we have after the first time step ($n \geq 1$):
\be
\bal
&\frac{\rho^{n+\fl 1 2}-\rho^n}{\Delta t/2} =\nabla_{\vx} \cdot  \left(\frac{1}{3\sigma^S} \nabla_{\vx}\rho^{n-\fl 1 2}\right)-\sigma^A\rho^{n}+G,\\
&\frac{\rho^{n+1} -\rho^n}{\Delta t}=\nabla_{\vx} \cdot  \left(\frac{1}{3\sigma^S} \nabla_{\vx}\rho^{n+\fl 1 2}\right) -  \sigma^A\rho^{n+\fl 1 2}+  G,
\eal
\ee
which is a second-order explicit RK scheme for the limiting diffusion equation (\ref{eqn:limiting1}). Therefore, the scheme is AP.

\begin{remark}
There are many other choices to construct the second order scheme by altering the order of $K$, $L$ and $S$, as long as the steps 2-4 are symmetric with respect to steps 5-7. Note that the AP property is always guaranteed due to the flexibility in the first order scheme.
\end{remark}

\subsection{Fully discrete scheme}

It remains for us to specify the discretization in the physical space and velocity space. This is the purpose of this section.

\subsubsection{Velocity discretization}

For the velocity space $\mathbb{S}^2$, we adopt the discrete velocity method\footnote{In the context of radiative transfer, this is usually referred to as discrete ordinates or $S_N$ method.}. The velocity points $\{\vv_i\}_{i=1,\dots,N_{\vv}}$ and weights $\{w_i\}_{i = 1,\ldots,N_{\vv}}$ are chosen according to the Lebedev quadrature on $\mathbb{S}^2$. Then all the integrals of the form $\langle F(\vv) \rangle_{\vv}$ are approximated as
\be
	\langle F(\vv)\rangle_{\vv} \approx \sum\limits_{i=1}^{N_{\vv}} w_i F(\vv_i).
\ee

\subsubsection{Spatial discretization}
\label{subsec:space}

For the physical space $\Omega_{\vx}$, we assume the third dimension is homogeneous and the domain is rectangular so that we consider $\vx=(x,y)\in [a,b]\times[c,d]$. For simplicity, we assume periodic boundary condition.

To obtain the asymptotic limit in a more compact stencil, we adopt the 2D staggered grid proposed in \cite{KFJ16}. We divide the $x$ and $y$ directions uniformly into $N_x$ and $N_y$ cells with size $\Delta x=(b-a)/N_x$, $\Delta y=(c-d)/N_y$, respectively. We denote the vertices by $x_k=a+k\Delta x$, $y_l=c+l\Delta y$ ($k=0,\dots, N_x$, $l=0,\dots,N_y$), and the cell centers by $x_{k+\frac 12}=a+(k+\frac 12)\Delta x$, $y_{l+\frac 12}=c+(l+\frac 12)\Delta y$ ($k=0,\dots, N_x-1$, $l=0,\dots,N_y-1$). We then place the unknowns $\rho$ and $g$ as in Figure~\ref{fig:StaggeredGrid}. Namely,
\begin{itemize}
\item  $\rho$ is located at the vertices $(x_k,y_l)$ and cell centers $(x_{k+\frac12},y_{l+\frac 1 2})$, i.e., the red dots in the figure;
\item  $g$ (hence $\{K_i, X_i\}_{i = 1,\ldots,r}$) is located at the face centers $(x_{k+\frac12},y_l)$ and $(x_k,y_{l+\frac12})$, i.e., the blue diamonds in the figure.
\end{itemize}

\usetikzlibrary{shapes,arrows}

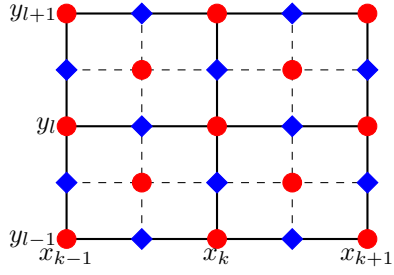
\begin{figure}[htb]
\centering
	\begin{tikzpicture} [x=2cm,y=1.5cm]
    \foreach \i in {0,...,2} {
        \draw [black, thick] (\i,0) -- (\i,2);
        \draw [black, thick] (0,\i) -- (2,\i) ;
    }
    \foreach \i in {1,...,2} {
        \draw [black, dashed] (\i-1/2,0) -- (\i-1/2,2) ;
        \draw [black, dashed] (0,\i-1/2) -- (2,\i-1/2);
    }
   \foreach \x in {1,...,2}{
      \foreach \y in {1,...,2}{
        \node[draw,circle,inner sep=2.5pt,fill,red] at (\x-1/2,\y-1/2) {};
      }
    }
    \foreach \x in {0,...,2}{
      \foreach \y in {1,...,2}{
        \node[draw,diamond,inner sep = 2.pt,fill,blue] at (\x,\y-1/2) {};
      }
    }
    \foreach \x in {1,...,2}{
      \foreach \y in {0,...,2}{
        \node[draw,diamond,inner sep = 2.pt,fill,blue] at (\x-1/2,\y) {};
      }
    }
        \foreach \x in {0,...,2}{
      \foreach \y in {0,...,2}{
        \node[draw,circle,inner sep=2.5pt,fill,red] at (\x,\y) {};
      }
    }
    \node [below] at (0,0) {$x_{k-1}$};
    \node [below] at (1,0) {$x_{k}$};
    \node [below] at (2,0) {$x_{k+1}$};
    \node [left] at (0,0) {$y_{l-1}$};
    \node [left] at (0,1) {$y_{l}$};
    \node [left] at (0,2) {$y_{l+1}$};
	\end{tikzpicture}
\caption{The staggered grids. $\rho$ is located at the red dots; $g$ (hence $\{K_i, X_i\}_{i = 1,\ldots,r}$) is located at the blue diamonds.} \label{fig:StaggeredGrid}
\end{figure}

In the following, we describe a second order finite difference method in space. We use simplified notations such as $\rho_{k,l} = \rho(x_k,y_l)$, $(K_i)_{k+\frac 12,l} = K_i(x_{k+\frac 12},y_l)$ to denote numerical solutions evaluated at the corresponding grid points. We use the first order $K$-$L$-$S$-$\rho$ scheme in time. The discussion for other time discretization methods is similar.

\begin{itemize}
\item $K$-step

Note that the system (\ref{semistepK}), in matrix form, can be written as
\be
\label{semistep1kreducedmatrixform}
	\fl {\mK^{n+1}-\mK^{n}}{\dt t} =-\mV^1 \p_{x} \mK^n -\mV^2 \p_{y} \mK^n+\dots,
\ee
where $\mK^n = [K_1^n,K_2^n,\ldots,K_r^n]^T$ and
\be
\mV^1_{jl} =   -\fl 1 \varepsilon  \left(\langle \xi V_j^{n}V_l^{n}\rangle_{\vv}-\fl 1 {4\pi}\langle V_j^{n}\rangle_{\vv}\langle \xi V_l^{n}\rangle_{\vv}\right), \quad \mV^2_{jl} =   -\fl 1 \varepsilon  \left(\langle \eta V_j^{n}V_l^{n}\rangle_{\vv}-\fl 1 {4\pi}\langle V_j^{n}\rangle_{\vv}\langle \eta V_l^{n}\rangle_{\vv}\right).
\ee
It is clear that the matrices $\mV^1$ and $\mV^2$ are not necessarily symmetric hence the system might not be hyperbolic. Therefore, to get a reasonable spatial discretization for (\ref{semistepK}), we propose to discretize the original equation (\ref{eqn:RHS}) and then project the resulting scheme. 

Specifically, we first discretize (\ref{eqn:RHS}) as
\be
\label{gsemi}
\bal
\p_tg=&-\fl 1 \varepsilon\left(I-\fl 1 {4\pi}\langle \ \rangle_{\vv}\right)\left(\xi^{+}  D_{+}^x g+\xi^{-} D_{-}^x g \right)-\fl 1 \varepsilon\left(I-\fl 1 {4\pi}\langle \ \rangle_{\vv}\right)\left(\eta^{+}D_{+}^y g+\eta^{-} D_{-}^y g \right)\\
&-\fl 1 {\varepsilon^2}\left(\xi D_c^x \rho+\eta D_c^y \rho\right)-\left(\fl {\sigma^S} {\varepsilon^2} +\sigma^A\right)g, 
\eal
\ee
where $\xi^+ = \text{max}(0,\xi)$, $\xi^- = \text{min}(0,\xi)$. A second order upwind operator is applied to the spatial derivatives of $g$ and a central difference operator is applied to the spatial derivatives of $\rho$. More precisely, we use
\be
\bal
D_+^x g(x,y) &= \fl {3g(x,y) -4g(x-\dt x,y) +g(x-2\dt x,y) }{2 \dt x}, \\
D_-^x g(x,y) &= \fl {-3g(x,y) +4g(x+\dt x,y) -g(x+2\dt x,y) }{2 \dt x},
\eal
\ee
and 
\be
D_c^x\rho(x,y) = \fl {\rho(x+\fl 1 2 \dt x,y)-\rho(x-\fl 1 2 \dt x,y)}{\dt x}.
\ee
Derivatives in $y$ are defined similarly.

We then project the equation (\ref{gsemi}) onto the space spanned by $\{V_j\}_{1\leq j\leq r}$, which yields
\be
\label{semistep1kreducedprojected}
\bal
	\fl {(K_j^{n+1})_{k+\fl 1 2,l}-(K_j^{n})_{k+\fl 1 2,l}}{\dt t} =& -\fl 1 \varepsilon \sum_{i=1}^r \left(\langle \xi^+ V_j^{n}V_i^{n}\rangle_{\vv}-\fl 1 {4\pi}\langle V_j^{n}\rangle_{\vv}\langle \xi^+ V_i^{n}\rangle_{\vv}\right) D_+^x (K_i^n)_{k+\fl 1 2,l} \\
	&-\fl 1 \varepsilon \sum_{i=1}^r \left(\langle \xi^- V_j^{n}V_i^{n}\rangle_{\vv}-\fl 1 {4\pi}\langle V_j^{n}\rangle_{\vv}\langle \xi^- V_i^{n}\rangle_{\vv}\right) D_-^x (K_i^n)_{k+\fl 1 2,l} \\
	&-\fl 1 \varepsilon \sum_{i=1}^r \left(\langle \eta^+ V_j^{n}V_i^{n}\rangle_{\vv}-\fl 1 {4\pi}\langle V_j^{n}\rangle_{\vv}\langle \eta^+ V_i^{n}\rangle_{\vv}\right) D_+^y (K_i^n)_{k+\fl 1 2,l} \\
	&-\fl 1 \varepsilon \sum_{i=1}^r \left(\langle \eta^- V_j^{n}V_i^{n}\rangle_{\vv}-\fl 1 {4\pi}\langle V_j^{n}\rangle_{\vv}\langle \eta^- V_i^{n}\rangle_{\vv}\right) D_-^y (K_i^n)_{k+\fl 1 2,l} \\
	&-\fl 1 {\varepsilon^2} \langle\xi V_j^n\rangle_{\vv}D_c^x\rho^{n}_{k+\fl 1 2,l}-\fl 1 {\varepsilon^2} \langle\eta V_j^n\rangle_{\vv}D_c^y\rho^{n}_{k+\fl 1 2,l}\\
	&-\fl {\sigma^S_{k+\fl 12,l}}{\varepsilon^2}(K_j^{n+1})_{k+\fl 1 2,l}-\sigma^A_{k+\fl 12,l}(K_j^{n})_{k+\fl 1 2,l}.
\eal
\ee
Here the scheme is given at the grid points $(x_{k+\fl 1 2},y_l)$. The scheme at the grid points $(x_k,y_{l+\frac12})$ is similar.

\item $L$-step and $S$-step

One can add spatial discretization to (\ref{semistepL}) and (\ref{semistepS}) directly. First of all, we approximate the inner product $\langle  \ \rangle_{\vx}$ by a midpoint rule:
\be
	\langle F(x,y)\rangle_{\vx} =  \int_{{[a,b]^2}}F  \, {\rm d} x {\rm d} y \approx
	\frac{1}{2}\Delta x \Delta y \sum\limits_{k=1}^{N_x}\sum\limits_{l=1}^{N_y} (F_{k+\frac12,l}+F_{k,l+\frac12}).
\ee
Then we approximate the spatial derivatives of $\rho$ and $X_i$ at $(x_{k+\fl 1 2},y_l)$ and $(x_k,y_{l+\frac12})$ by
\begin{align}
&\p_x\rho_{k+\fl 1 2,l} \approx \fl {\rho_{k+1,l}-\rho_{k,l}}{\dt x}, \quad  &&\p_x (X_{i})_{k+\fl 1 2,l} \approx \fl {(X_i)_{k+\fl 3 2,l} -(X_i)_{k-\fl 1 2,l}}{2 \dt x}, 	\\
&\p_x\rho_{k,l+\fl 1 2} \approx \fl {\rho_{k+\fl 1 2,l+\fl 1 2}-\rho_{k-\fl 1 2,l+\fl 1 2}}{\dt x},\quad &&\p_x (X_i)_{k, l+\fl 1 2} \approx \fl { (X_i)_{k+1,l+\fl 1 2} -(X_i)_{k-1,l+\fl 1 2}}{2 \dt x}.
\end{align}
Derivatives in $y$ are treated similarly.

\item $\rho$-step

At the grid points $(x_k,y_l)$, (\ref{semisteprho}) is discretized as
\be \label{discreterho}
\bal
\fl {\rho^{n+1}_{k,l}-\rho^{n}_{k,l}}{\dt t} =&-\fl 1 {4\pi}\sum\limits_{i,j=1}^r \fl {(X_i^{n+1})_{k+\fl 1 2,l} -(X_i^{n+1})_{k-\fl 1 2,l}}{ \dt x}  S_{ij}^{n+1} \langle\xi V_j^{n+1}\rangle_{\vv} \\
& -\fl 1 {4\pi}\sum\limits_{i,j=1}^r \fl {(X_i^{n+1})_{k,l+ \fl 12} -(X_i^{n+1})_{k,l-\fl 12}}{ \dt y}  S_{ij}^{n+1} \langle\eta V_j^{n+1}\rangle_{\vv}-\sigma^A_{k,l}\rho^{n}_{k,l}+G_{k,l}.
\eal
\ee
The scheme at the grid points $(x_{k+\fl 12},y_{l+\frac12})$ is similar.
\end{itemize}

\subsubsection{AP property of the fully discrete scheme}

Similar to the semi-discrete case, in the limit $\varepsilon \rightarrow 0$, the $K$-$L$-$S$ steps yield
\be
\bal
	&g^{n+1}_{k+\fl 1 2,l}=\sum_{i,j=1}^r (X_i^{n+1})_{k+\fl 12,l}S_{ij}^{n+1}V_j^{n+1}=-\fl 1 {\sigma^S_{k+\fl 1 2,l}}\left(\xi  \fl {\rho^n_{k+1,l}-\rho^n_{k,l}}{\dt x}+\eta \fl{\rho^{n}_{k+\fl 1 2,l+\fl 1 2}-\rho^n_{k+\fl 1 2,l-\fl 1 2}}{\dt y}\right),\\
	&g^{n+1}_{k,l+\fl 1 2}=\sum_{i,j=1}^r (X_i^{n+1})_{k,l+\fl 12}S_{ij}^{n+1}V_j^{n+1}=-\fl 1 {\sigma^S_{k,l+\fl 1 2}}\left(\xi  \fl {\rho^n_{k+\fl 1 2,l+\fl 1 2}-\rho^n_{k-\fl 1 2,l+\fl 1 2}}{\dt x}+\eta \fl{\rho^n_{k,l+1}-\rho_{k,l}^n}{\dt y}\right),
\eal
\ee
which, when substituting into (\ref{discreterho}), give
\be
\bal
\label{diffusionlimit1}
\fl {\rho^{n+1}_{k,l}-\rho^{n}_{k,l}}{\dt t}=&\fl 1 3 \fl {1}{\dt x^2}\left(\fl {\rho^{n}_{k+1,l}-\rho^{n}_{k,l}} { {\sigma^S_{k+\fl 1 2,l}}}-\fl {\rho^{n}_{k,l}-\rho^{n}_{k-1,l}} { {\sigma^S_{k-\fl 1 2,l}}}\right)+\fl 1 3 \fl {1}{\dt y^2}\left(\fl {\rho^{n}_{k,l+1}-\rho^{n}_{k,l}} { {\sigma^S_{k,l+\fl 1 2}}}-\fl {\rho^{n}_{k,l}-\rho^n_{k,l-1}} { {\sigma^S_{k,l-\fl 1 2}}}\right)\\
&-\sigma^A_{k,l}\rho^{n}_{k,l}+G_{k,l}.
\eal
\ee
This is an explicit standard 5-point finite difference scheme applied to the limiting diffusion equation (\ref{eqn:limiting1}) at grid points $(x_k,y_l)$. The limiting scheme at grid points $(x_{k+\fl 1 2},y_{l+\fl 1 2})$ can be considered similarly. Therefore, the fully discrete scheme is also AP.

\section{A Fourier analysis of the low-rank structure of the solution}
\label{sec:Fourier}

In this section, we analyze the behavior of the solution to the linear transport equation by performing a simple Fourier analysis. Our focus is in the kinetic regime because the rank is already proved to be small in the diffusive regime. 

For simplicity, we consider the 1D slab geometry $x\in[0,2\pi]$ with periodic boundary condition, and $v\in [-1,1]$ (so $\langle \ \rangle_v =\int_{-1}^1 \cdot \, \rd{v}$). Also we assume $\sigma^A=G=0$. Then the macro-micro system of the linear transport equation reads:
\be
\bal
&\p_t\rho =-\fl 1 {2}\langle v  \p_x g\rangle_{v},\\
&\p_tg=-\fl 1 \varepsilon\left(I-\fl 1 {2}\langle \ \rangle_{v}\right) (v \p_x g)-\fl 1 {\varepsilon^2}v \p_x \rho-\fl {\sigma^S} {\varepsilon^2} g.
\eal
\ee
Projecting the above system onto the Fourier space of $x$ yields 
\be \label{Fourier}
	\bal
	&\p_t \hat{\rho}_{k} =- \fl 1 {2} ik \langle v \hat{g}_{k}\rangle_{v},\\
	&\p_t\hat{g}_{k}=- \fl 1 {\varepsilon}ik\left(v\hat{g}_{k}-\fl 1 {2}\langle v \hat{g}_{k}\rangle_{v}\right)-\fl 1 {\varepsilon^2}ivk \hat{\rho}_{k}- \fl 1 {\varepsilon^2}\sum\limits_{k_1 = -\infty}^{\infty}\hat{g}_{k-k_1}
	\hat{\sigma}_{k_1},
	\eal
\ee
where $\hat{\rho}_k(t)$, $\hat{g}_k(t,v)$ and $\hat{\sigma}_{k}$ are the Fourier coefficients of $\rho$, $g$ and $\sigma^S$, respectively.

For a constant $\sigma^S$ we have
\be
\hat{\sigma}_0=\sigma^S, \quad \hat{\sigma}_k=0, \ k\neq 0,
\ee
and the system (\ref{Fourier}) reduces to
 \be
	\bal
	&\p_t \hat{\rho}_{k} =- \fl 1 {2} ik \langle v \hat{g}_{k}\rangle_{v},\\
	&\p_t\hat{g}_{k}=- \fl 1 {\varepsilon}ik\left(v\hat{g}_{k}-\fl 1 {2}\langle v \hat{g}_{k}\rangle_{v}\right)-\fl 1 {\varepsilon^2}ivk \hat{\rho}_{k}- \fl 1 {\varepsilon^2}\sigma^S\hat{g}_k.
	\eal
\ee
Hence all the frequency modes are decoupled. It is clear that if initially 
\be \label{bandinitial}
\rho(0,x)=\sum_{k=-m_0}^{m_0}\hat{\rho}_k(0) e^{ikx}, \quad g(0,x,v)=\sum_{k=-m_0}^{m_0}\hat{g}_k(0,v)e^{ikx},
\ee
i.e., $\rho(0,x)$ and $g(0,x,v)$ are band-limited, then the latter solutions will remain in the same frequency range. In this case the solution is clearly low-rank. 

This analysis is similar to the analysis conducted in \cite{piazzola}, where it was shown that for the  linearized Vlasov--Maxwell equation the solution remains low rank if it is initially in a form similar to \eqref{bandinitial}. However, the present situation is different in the sense that if we have a non-constant $\sigma^S$, as is commonly the case in practice, then even if the initial value is in that form additional Fourier modes are excited gradually with time. This is as far as we can go with such an argument.

However, it should not be taken to imply that performing a low-rank approximation is necessarily futile in such a situation. In fact, the dynamical low-rank integrator makes no assumptions that the space dependence has to take the form of a finite number of Fourier modes; this is purely an artifact of the analysis done here. Hence, just because we have an infinite number of Fourier modes does not necessarily imply the solution can not be captured by a low-rank scheme. In fact, from the numerical tests in the next section, we can see that the rank of the solution in the kinetic regime when $\sigma^S$ is not constant can be rather intricate, but that often relatively small ranks are sufficient to obtain an accurate approximation to the dynamics of interest. 

\section{Numerical results}
\label{sec:num}

In this section, we present several numerical examples to demonstrate the accuracy and efficiency of the proposed low rank method. In all examples, we consider a two-dimensional square domain in physical space, i.e.~$\vx = (x,y)\in [a,b]^2$ and periodic boundary conditions. Note that in some of the examples (e.g., the line source problem), one has to choose a large number of points in the angular direction to obtain a reasonable solution (for both the full tensor and low rank methods). This is the well-known drawback of the discrete velocity or collocation method. If a Galerkin rather than collocation approach is adopted, one could potentially use fewer discretization points (or bases). As the focus of the paper is on the low rank method, we leave the comparison of different angular discretizations to a future study. 
\subsection{Accuracy test}
\label{subsec:accuracy}

We first examine the accuracy of our method (in time and space) using a manufactured solution. We choose
\be
f(t,x,y,\xi,\eta,\gamma) = \exp(-t)\sin^2(2\pi x)\sin^2(2\pi y)\left(1+\varepsilon\left(\fl {\eta+\eta^3}{3}\right)\right), \quad (x,y) \in [0,1]^2.
\ee
The corresponding $\rho$ and $g$ are
\be
\bal
&\rho(t,x,y)=\exp(-t)\sin^2(2\pi x)\sin^2(2\pi y),\\
&g(t,x,y,\xi,\eta,\gamma)=\exp(-t)\sin^2(2\pi x)\sin^2(2\pi y)\left(\fl {\eta+\eta^3}{3}\right).
\eal
\ee 
Let the scattering and absorption coefficients be $\sigma^S = 1$, $\sigma^A = 0$, then the source term $G$ is given by
\begin{equation}
G(t,x,y,\xi,\eta,\gamma) = \partial_t f+\fl 1 {\varepsilon} \vv \cdot \nabla_{\vx} f +\fl {1}{\varepsilon}g.
\end{equation}
We use this source term and the initial condition $\rho(t=0,x,y)$ and $g(t=0,x,y,\xi,\eta,\gamma)$ as input for our low rank method and compute the solution up to a certain time. Note that the source term here depends also on time and velocity, hence the scheme needs to be modified accordingly to take into account this dependency. We omit the details.

We consider both the first order scheme in Section~\ref{subsec:first} and the second order scheme in Section~\ref{subsec:second}, coupled with the second order spatial discretization described in Section~\ref{subsec:space}. We always take $N_{\vv}=590$ Lebedev quadrature points on the sphere $\mathbb{S}^2$ \cite{lebedev}. Since we know {\it a priori} the rank of the exact solution $g$ is 1, we fix $r = 5$ in the low rank method which is certainly sufficient to obtain an accurate solution.

We vary the spatial size $\Delta x$ and the value of $\varepsilon$, and evaluate the error at $t=0.1$ as
\begin{equation}
\left(\Delta x^2\sum_{k,l=1}^{N_x} \left(\rho_{\text{low rank}}(x_{k+\frac{1}{2}},y_{l+\frac{1}{2}})-\rho_{\text{exact}}(x_{k+\frac{1}{2}},y_{l+\frac{1}{2}})\right)^2\right)^{\frac{1}{2}}.
\end{equation}
Since the proposed schemes are AP, we expect them to be stable under a hyperbolic CFL condition when $\varepsilon \sim O(1)$ and a parabolic CFL condition when $\varepsilon \ll 1$. Specifically, we consider three kinds of CFL conditions: mixed CFL condition \(\dt t \sim c_1 \dt x^2+c_2 \varepsilon\dt x\), hyperbolic CFL condition \(\dt t \sim  \dt x\), and parabolic CFL condition \(\dt t \sim  \dt x^2\).

\begin{figure}[!htp]
\begin{center}
\includegraphics[width=2.6in]{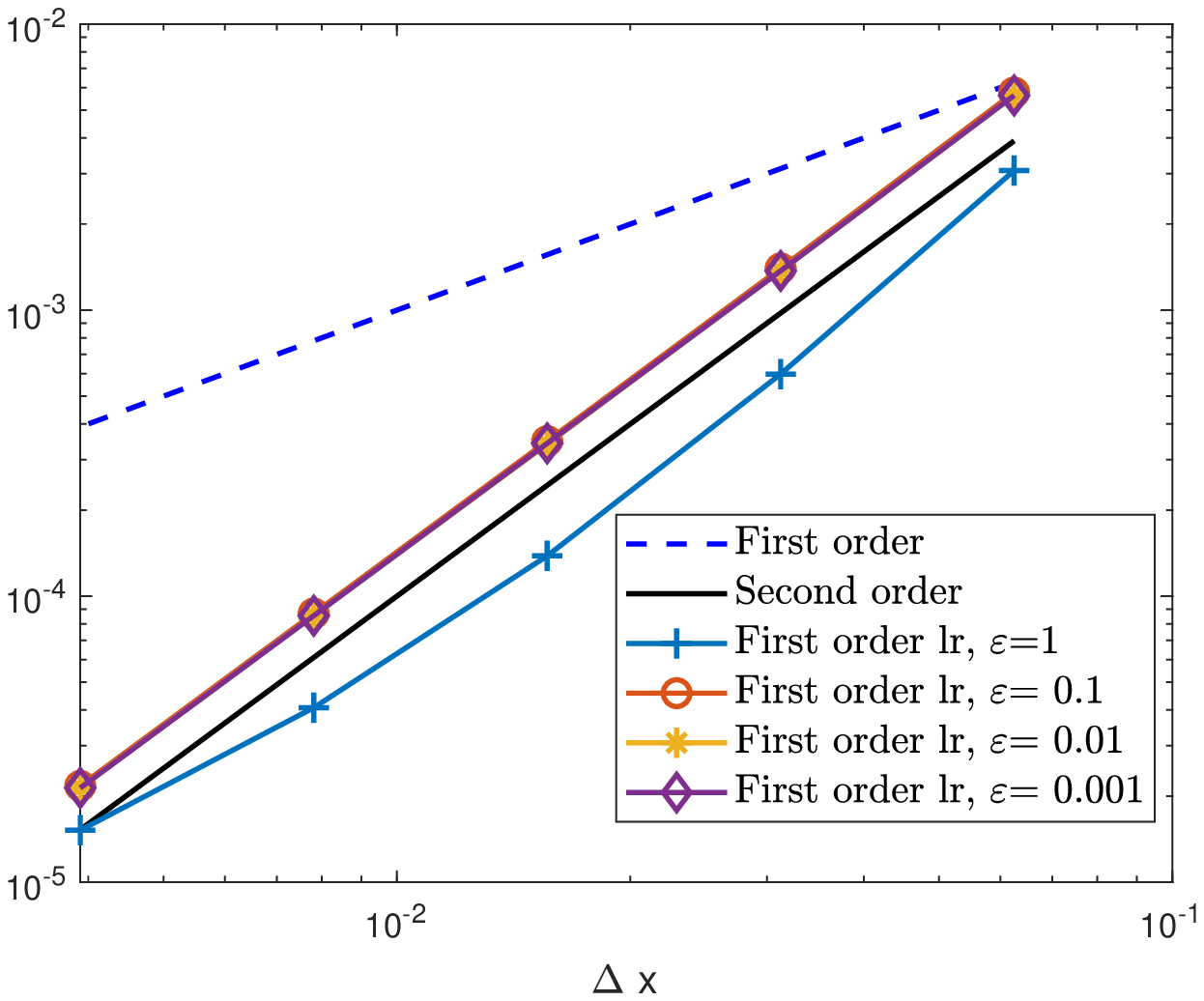}
\includegraphics[width=2.6in]{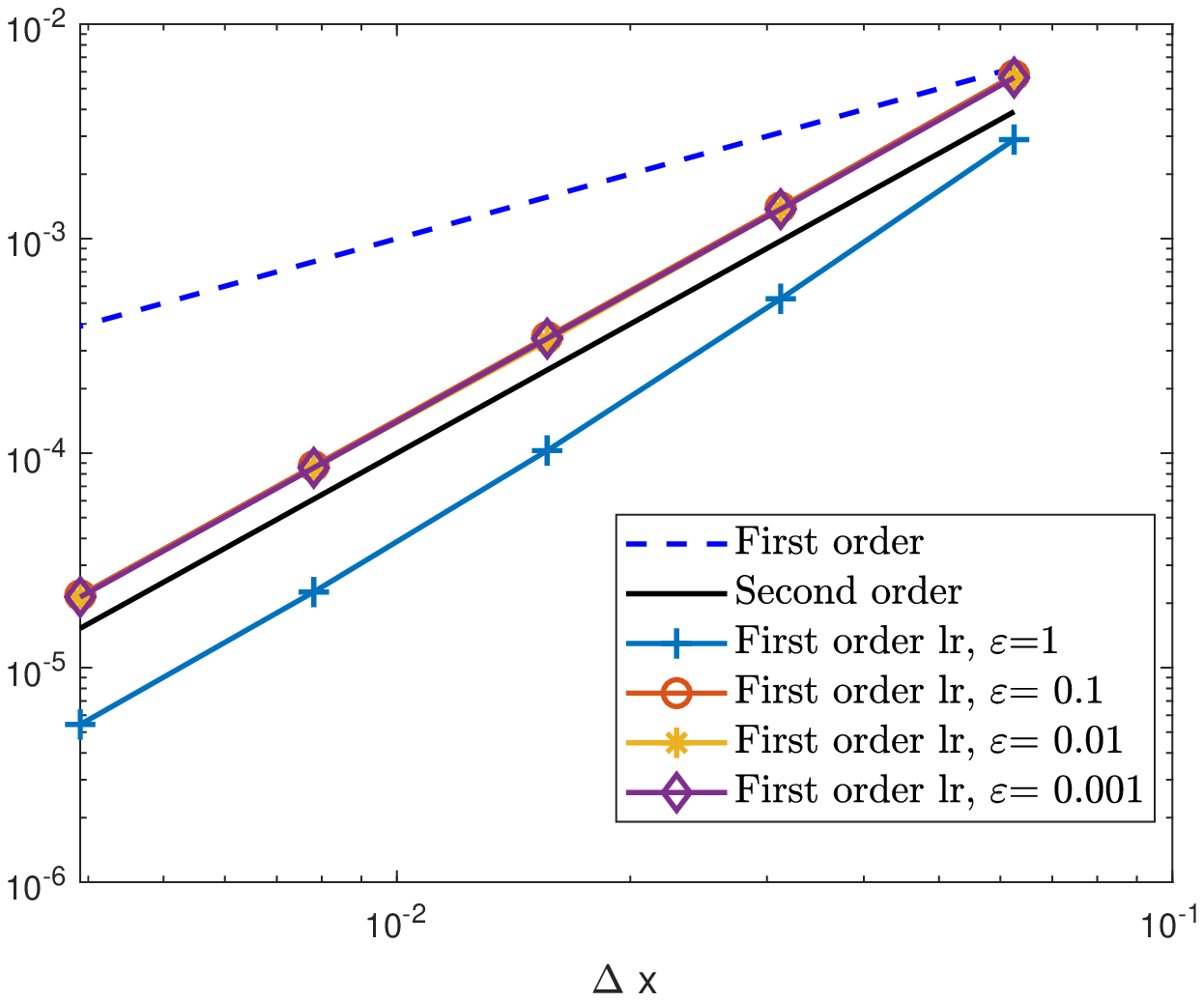}
\caption{Section~\ref{subsec:accuracy}: convergence order (first order low rank scheme). $l^2$-error v.s. $\dt x$. Left: mixed CFL condition $\dt t = 0.18\dt x^2+0.1\varepsilon\dt x$. Right: parabolic CFL condition $\dt t = 0.25\dt x^2$. Blue dashed line and black line are reference slopes of 1 and 2, respectively.}
\label{figure:1}
\end{center}
\end{figure}

The results of the first order (in time) scheme are shown in Figure~\ref{figure:1}. Under the mixed CFL condition, we expect to see first order convergence in the kinetic regime ($\varepsilon \sim O(1)$) and second order in the diffusive regime ($\varepsilon \ll 1$), which is clearly observed in Figure~\ref{figure:1} (left). Under the parabolic CFL condition, we always expect second order convergence, which is also clear in Figure~\ref{figure:1} (right).

\begin{figure}[!htp]
\begin{center}
\includegraphics[width=2.6in]{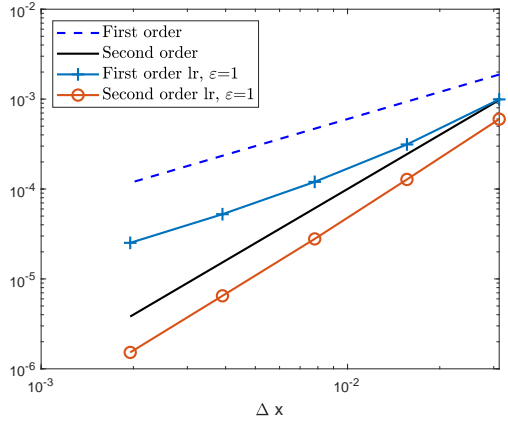}
\caption{Section~\ref{subsec:accuracy}: convergence order (second order low rank scheme). $l^2$-error v.s. $\dt x$. Hyperbolic CFL condition $\dt t = 0.4\dt x$ is used. Blue dashed line and black line are reference slopes of 1 and 2, respectively. Result of the first order scheme under the same CFL condition is plotted also for comparison.}
\label{figure:2}
\end{center}
\end{figure}

For the second order (in time) scheme, we don't expect order higher than two in the diffusive regime since $\Delta t\sim \Delta x^2$ and the error behaves as $O(\Delta t^2+\Delta x^2)=O(\Delta x^4+\Delta x^2)$. Hence we only test its performance in the kinetic regime ($\varepsilon \sim O(1)$) under the hyperbolic CFL condition, where $\Delta t \sim \Delta x$ and the error is $O(\Delta t^2+\Delta x^2)=O(\Delta x^2)$. The result is shown in Figure~\ref{figure:2}, where we can clearly see the uniform second order accuracy of the scheme (in contrast to the first order scheme).

\subsection{Test with Gaussian initial value}

In this test case, we consider a smooth Gaussian initial condition:
\be
f(t=0,x,y,\xi,\eta,\gamma) = \fl {1}{4\pi\varsigma^2}\exp\left(-\fl {x^2+y^2}{4\varsigma^2}\right), \quad \varsigma^2 = 10^{-2},\quad  (x,y) \in [-1,1]^2,
\label{eqn:initial_gauss}
\ee
with zero absorption coefficient and source term $\sigma^A = G=0$. 

\subsubsection{Constant scattering coefficient $\sigma^S$}
\label{subsec:const}

We first consider $\sigma^S \equiv 1$ and focus on the AP property of the proposed method. Therefore, we set $\varepsilon = 10^{-6}$ and compare our first order low rank method with the reference solution obtained by integrating (\ref{diffusionlimit1}), which solves the limiting diffusion equation directly. In the low rank method, we use $N_x=N_y=128$, $N_{\vv}=590$ Lebedev quadrature points on $\mathbb{S}^2$, and time step $\dt t = 0.1\dt x^2+0.1\varepsilon\dt x$, and fix the rank as $r=5$. In solving the diffusion equation, we use $N_x=N_y=512$ and time step $\dt t = 0.75 \dt x^2$. 

\begin{figure}[!htp]
\begin{center}
\includegraphics[width=2in]{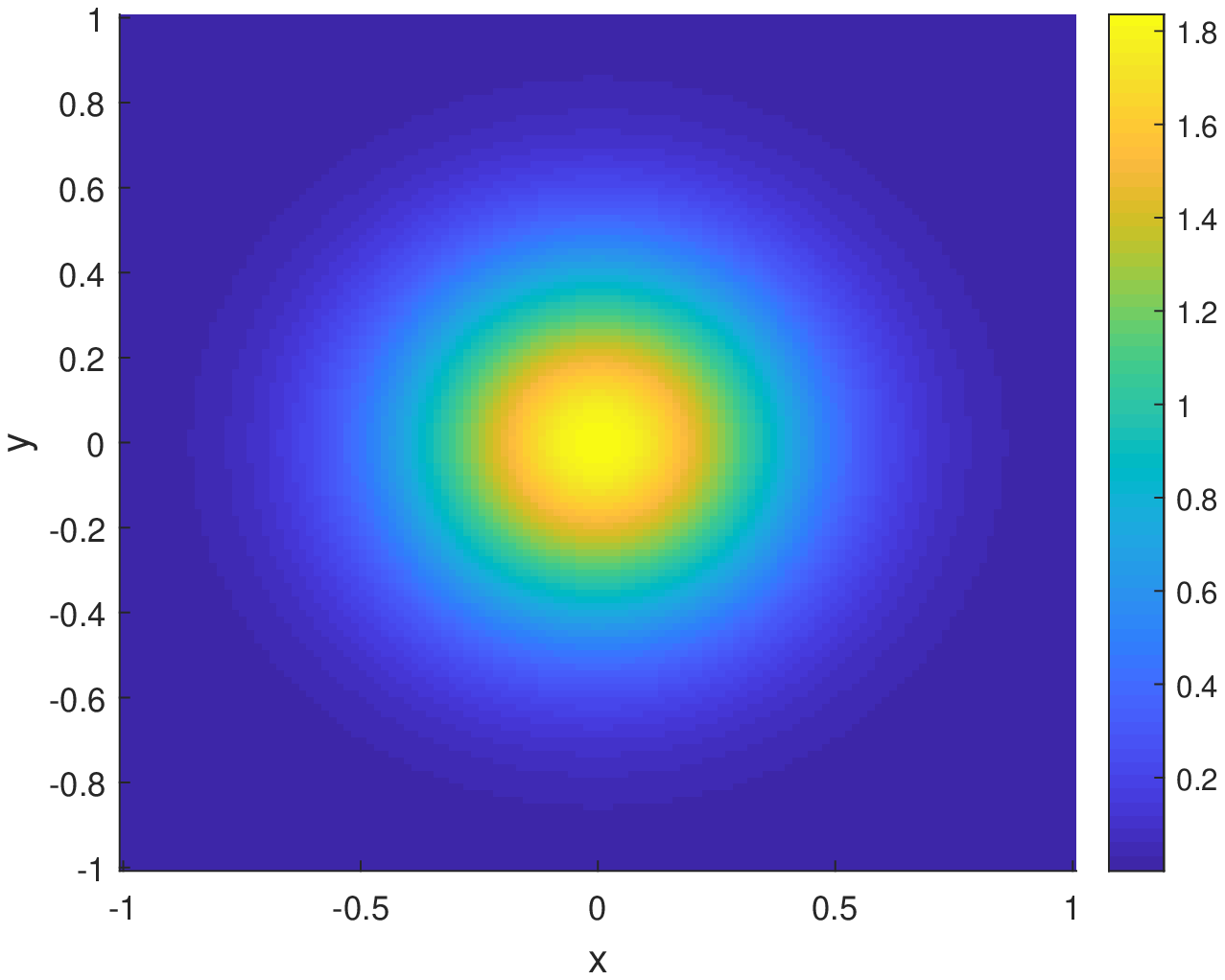}
\includegraphics[width=2in]{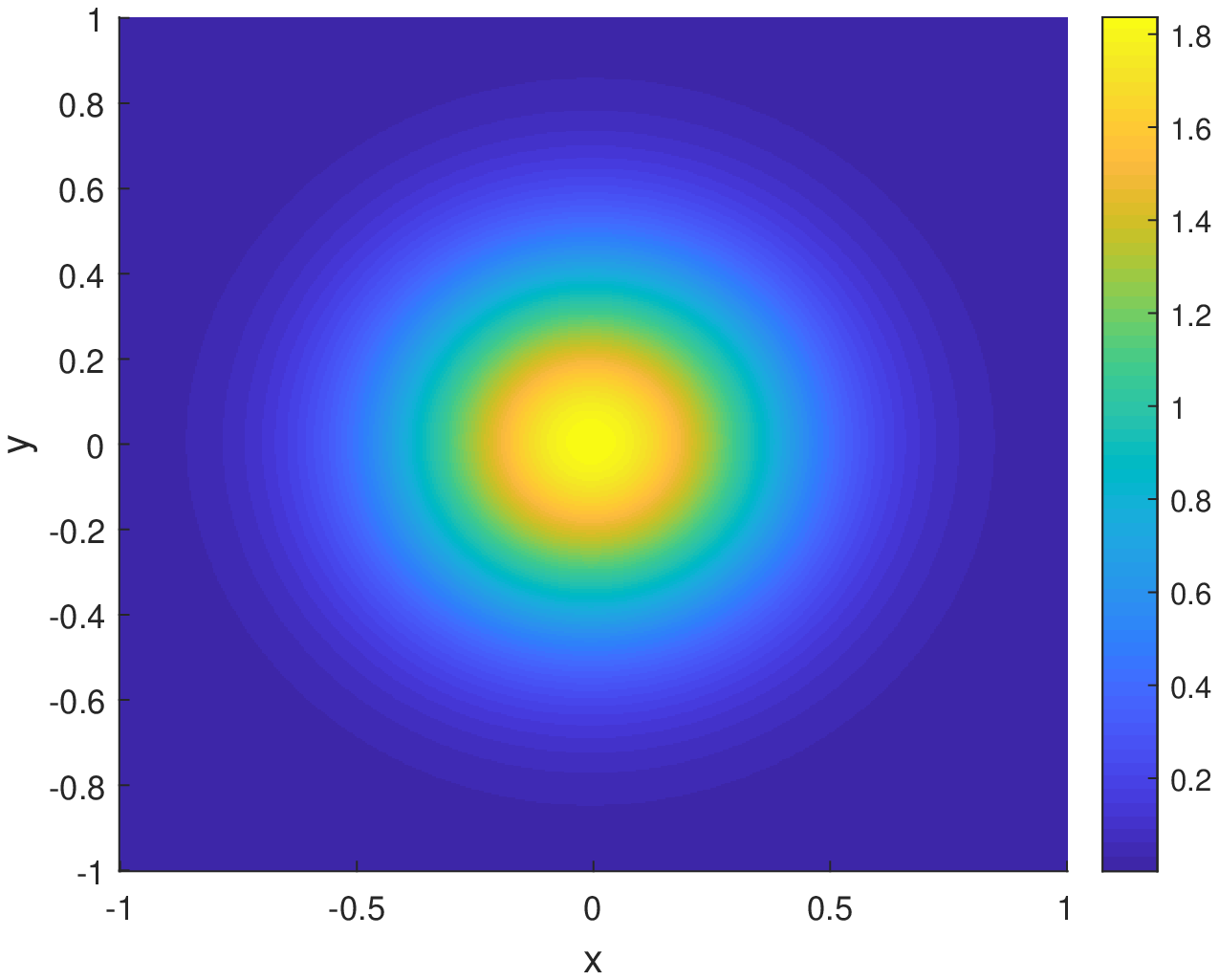}
\includegraphics[width=1.9in]{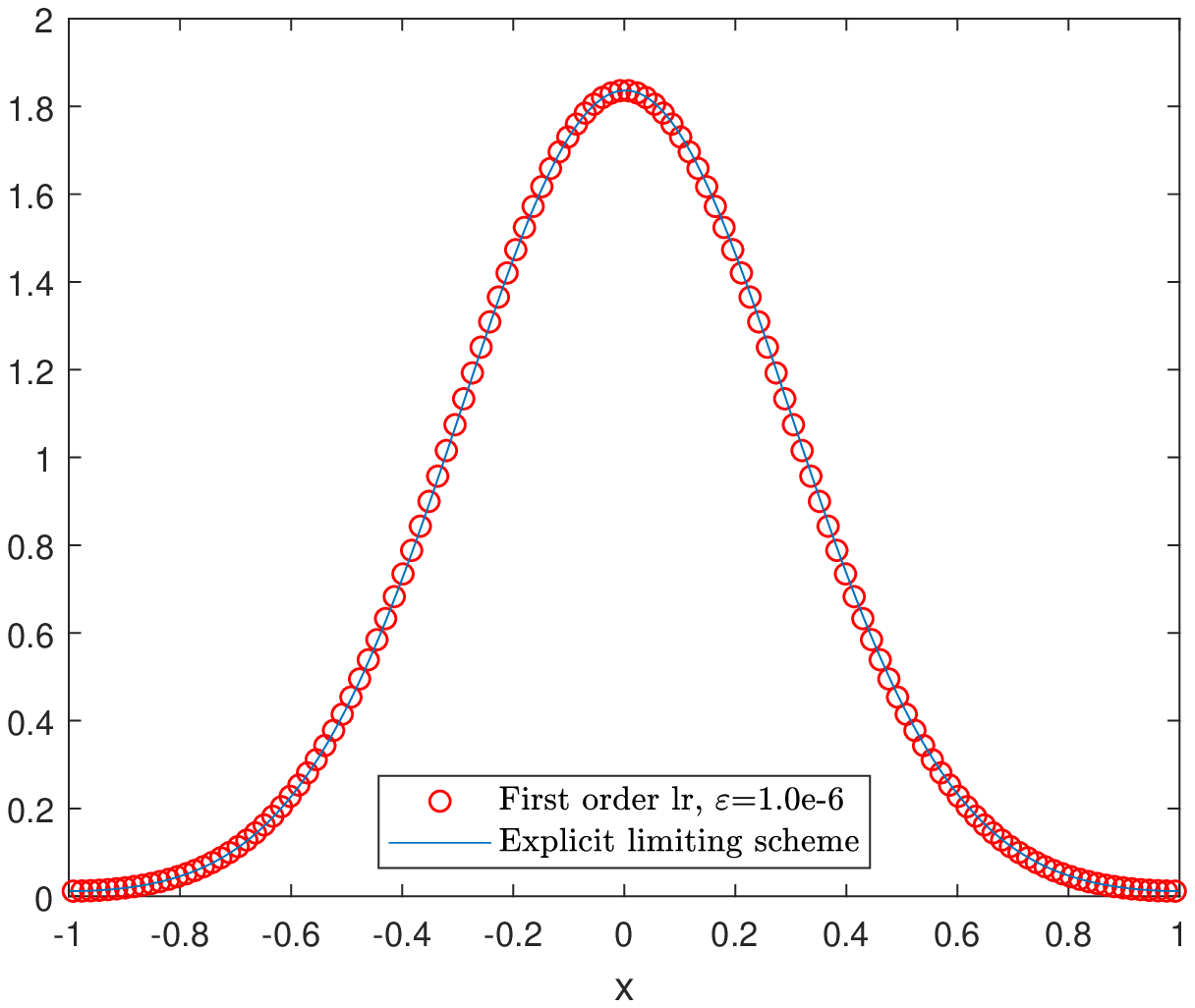}
\caption{Section~\ref{subsec:const}: constant scattering coefficient. Density profile of the low rank solution (left), reference solution to the limiting diffusion equation (middle), and comparison of two solutions with $y = 0$ (right).}
\label{figure:3}
\end{center}
\end{figure}

\begin{figure}[!htp]
\begin{center}
\includegraphics[width=2.6in]{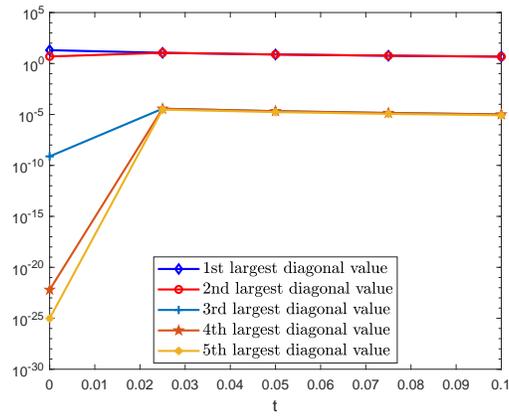}
\caption{Section~\ref{subsec:const}: constant scattering coefficient. Singular values of the matrix $S$ for the low rank method.}
\label{figure:3_1}
\end{center}
\end{figure}

The solutions at $t=0.1$ are shown in Figure~\ref{figure:3}, where they match very well. As the theory predicts, in the limiting diffusive regime, the solution $g$ should be become rank-2. To confirm this, we track the singular values of the matrix $S$, see Figure~\ref{figure:3_1}. Clearly, the effective rank is 2 (two singular values are above the threshold of $10^{-5}$, which is on the order of the spatial error).

\subsubsection{Variable scattering coefficient $\sigma^S$}
\label{subsec:variable}

\begin{figure}[!htp]
\begin{center}
\includegraphics[width=2.55in]{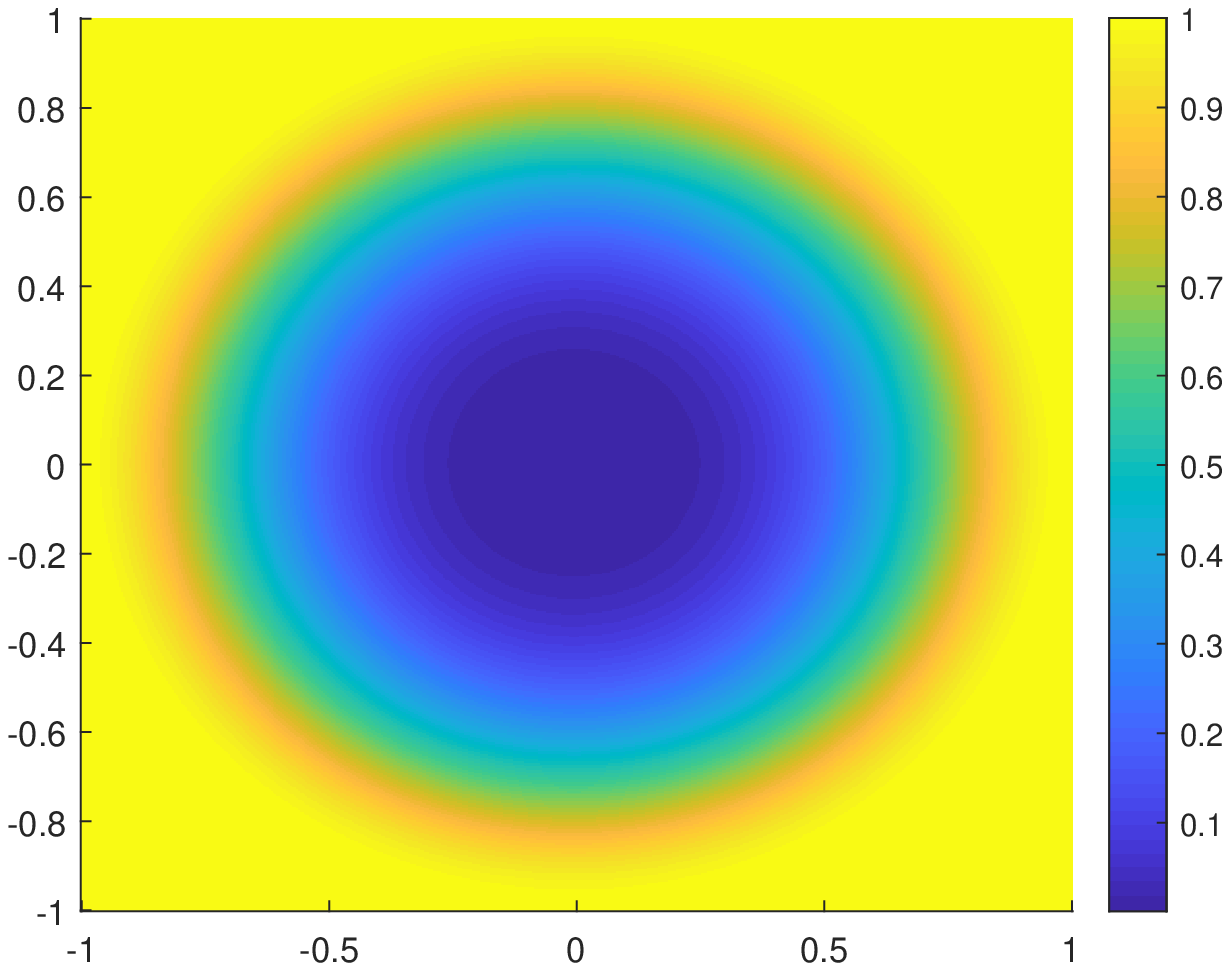}
\includegraphics[width=2.4in]{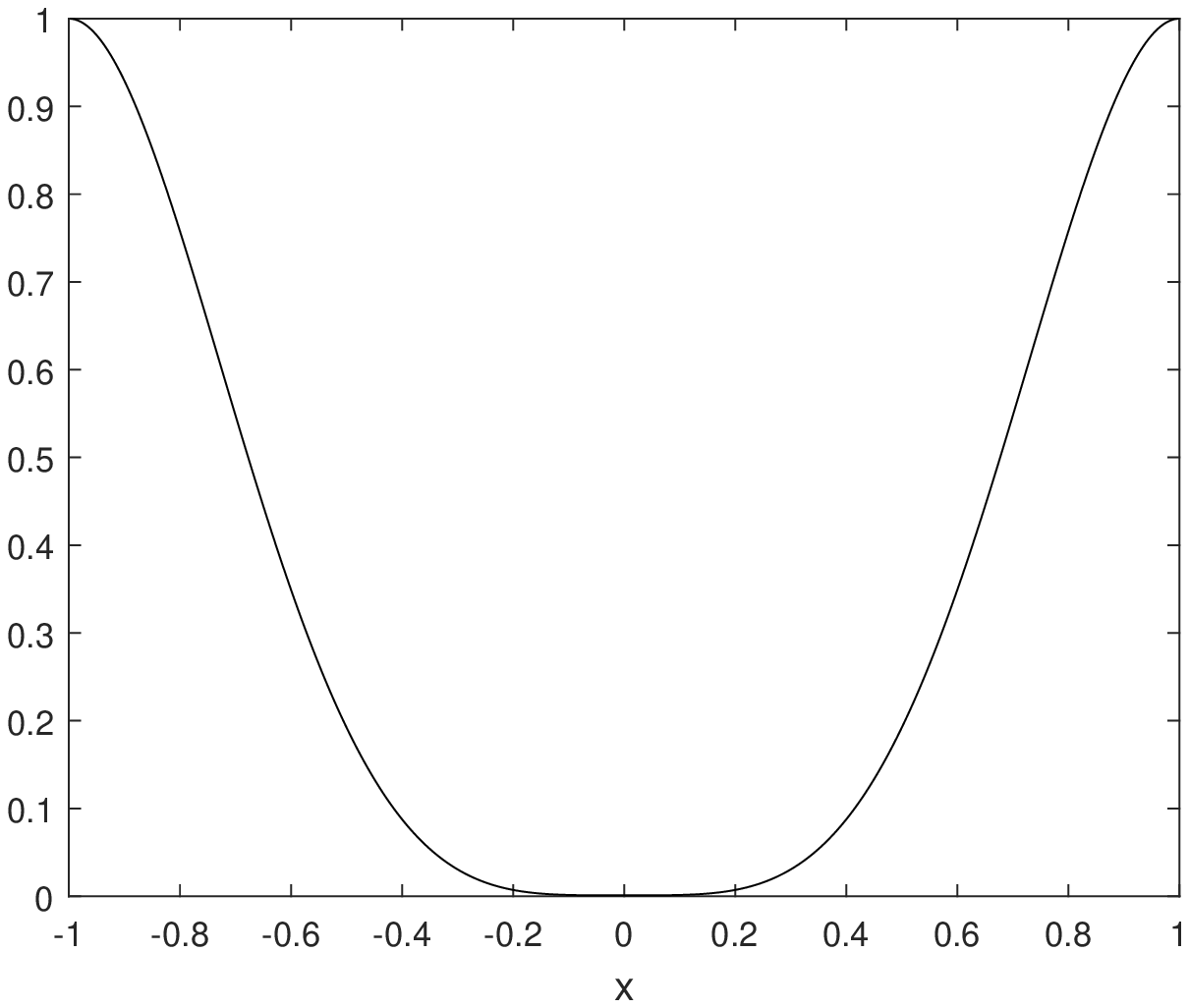}
\caption{Section~\ref{subsec:variable}: variable scattering coefficient. Profile of $\sigma^S$ (left) and a slice with $y = 0$ (right).}
\label{figure:4_0}
\end{center}
\end{figure}

We then set $\varepsilon = 0.01$ (an intermediate regime) and consider a spatially dependent scattering coefficient
\be
\sigma^S(x,y) = \left\{ \bal&0.999c^4(c+\sqrt{2})^2(c-\sqrt{2})^2+0.001,&\ &c=\sqrt{x^2+y^2}<1,&\\&1,&&\text{otherwise,}&\eal \right.
\ee
whose profile is shown in Figure~\ref{figure:4_0}. This is a challenging test as $\frac{\sigma^S(x,y)}{\varepsilon}$ varies in a large range $[0.1,100]$. Our aim here is to investigate the rank dependence of the low rank method and its performance compared with the full tensor method.

Specifically, we compare the first order low rank method with the first order IMEX method that solves the macro-micro decomposition of the linear transport equation directly \cite{LM08} (referred to as the full tensor method in the following). We use the same spatial mesh, same CFL condition $\dt t = 0.1\min(\sigma^S)\dt x^2+0.1\varepsilon\dt x$, and same $N_{\vv}=2702$ Lebedev quadrature points on $\mathbb{S}^2$ for both methods. In the low rank method, we choose different ranks from $20$ to $120$. 

\begin{figure}[!htp]
\begin{center}
\includegraphics[width=1.9in]{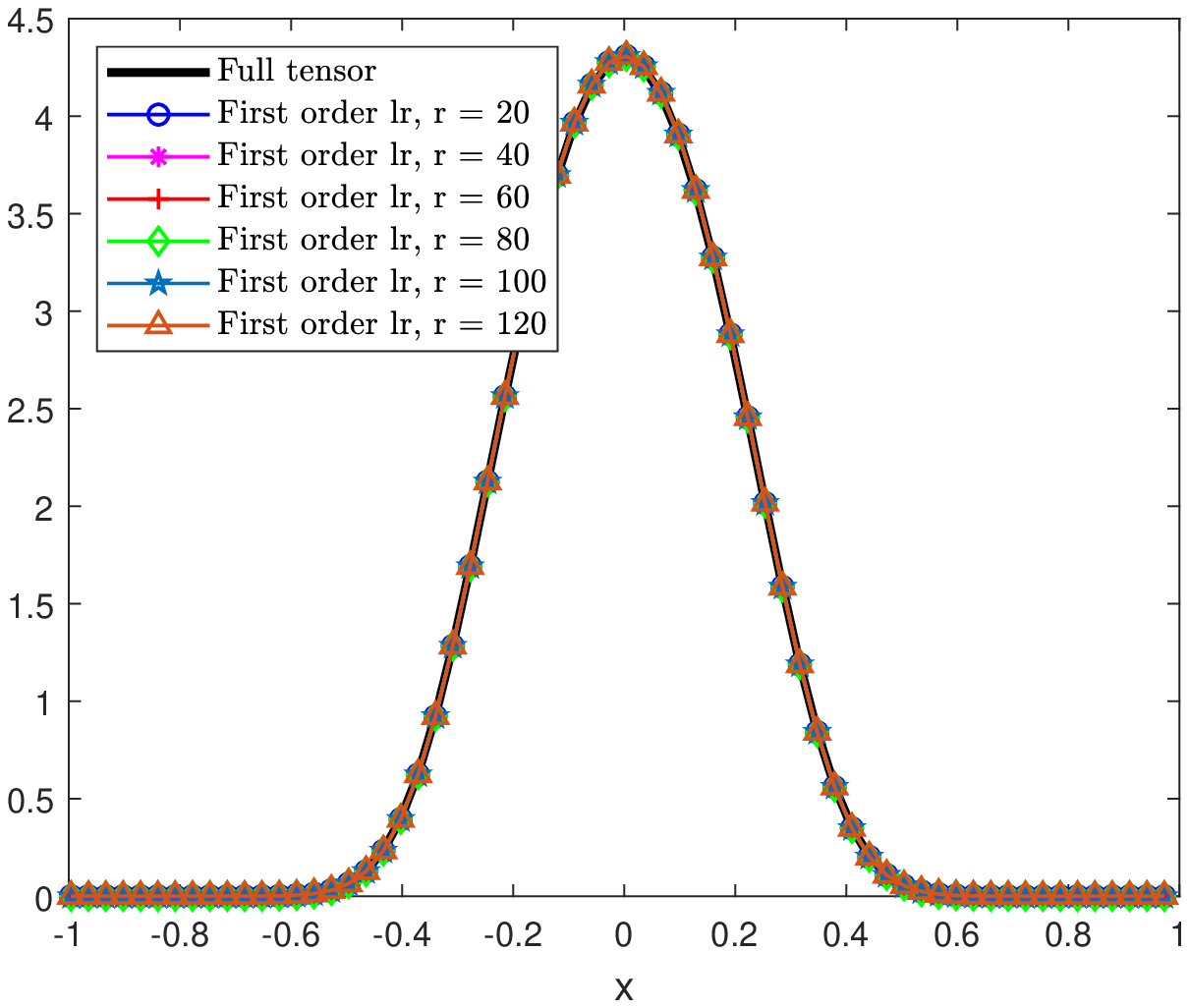}
\includegraphics[width=1.9in]{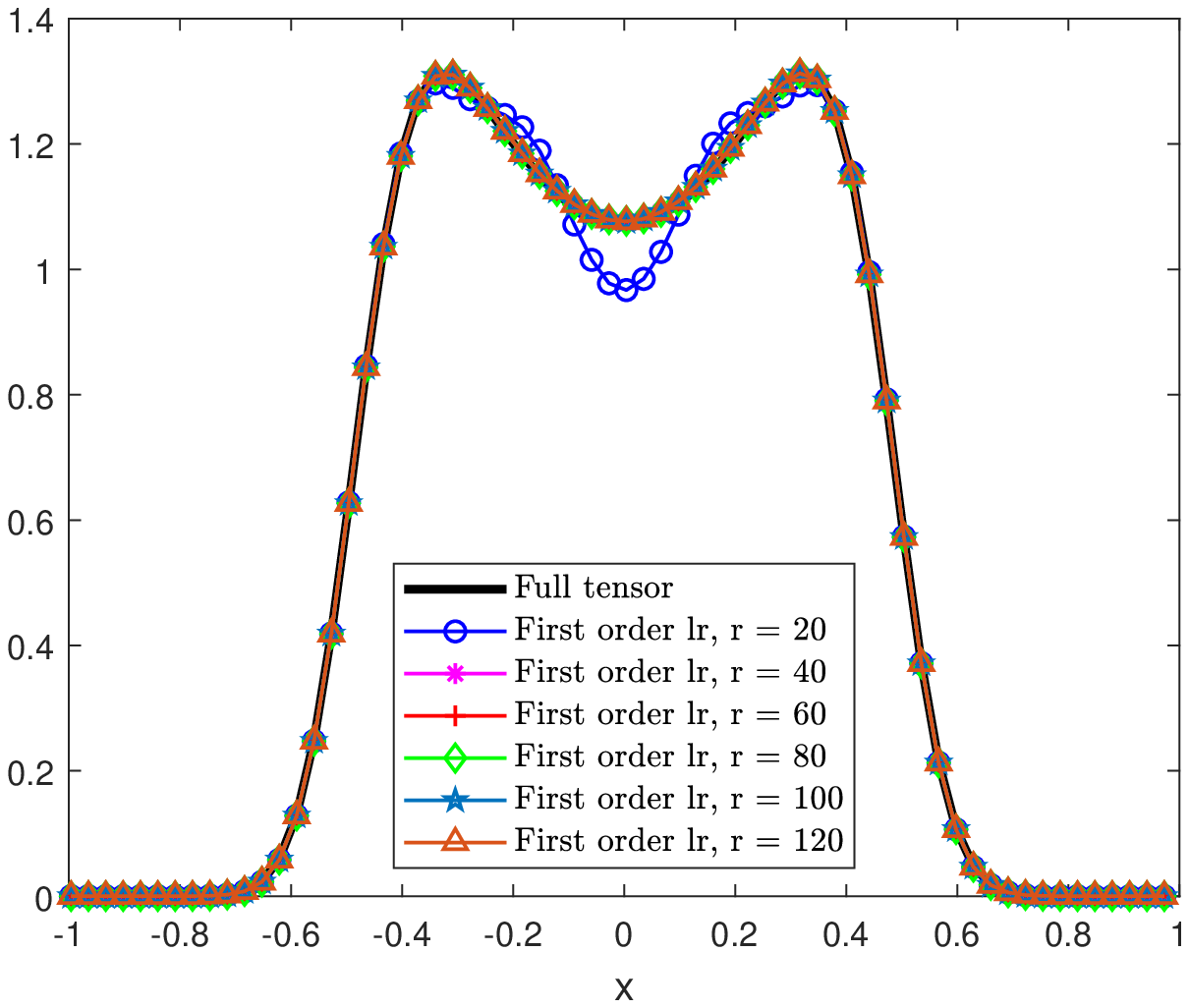}
\includegraphics[width=1.9in]{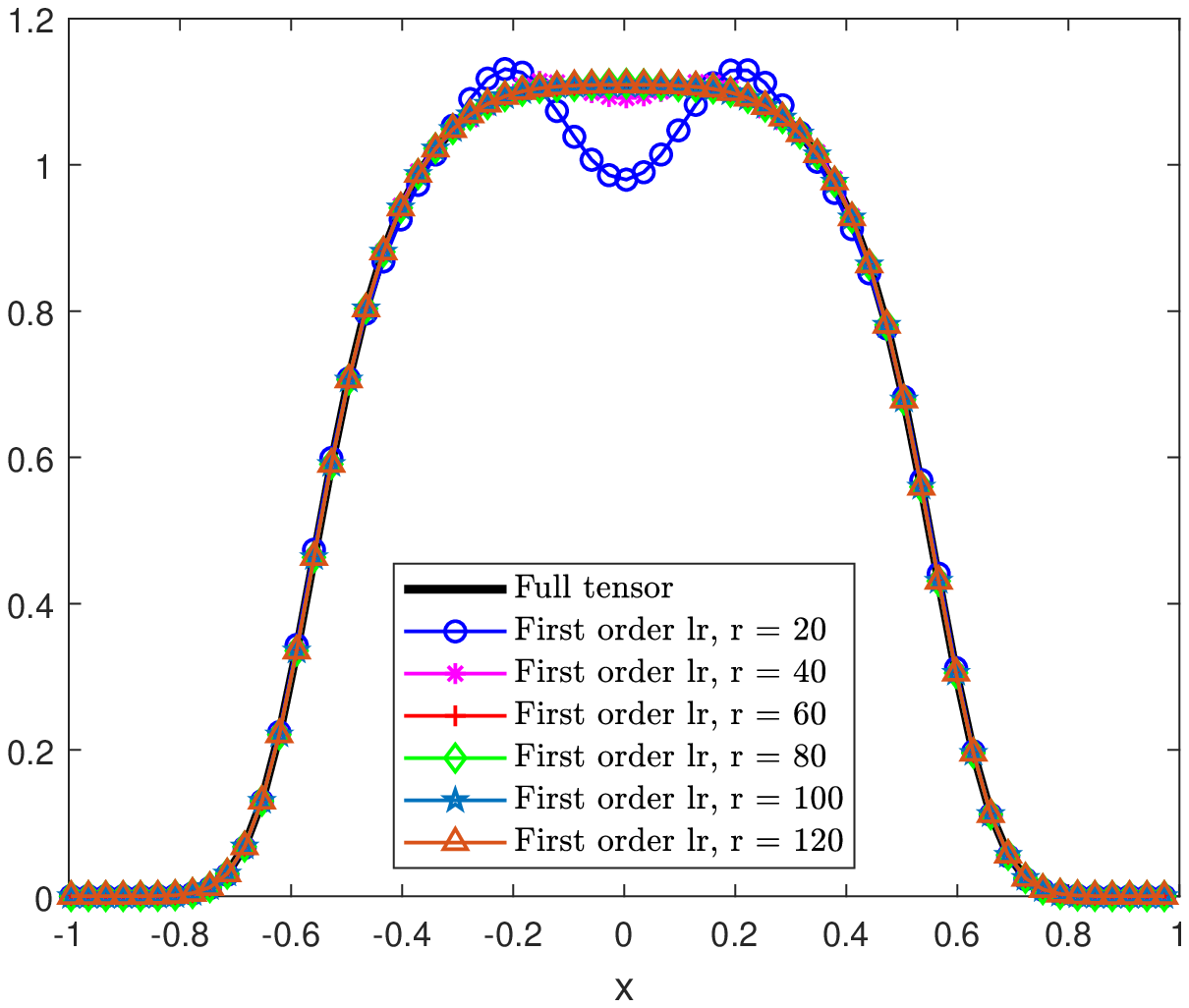}
\includegraphics[width=1.9in]{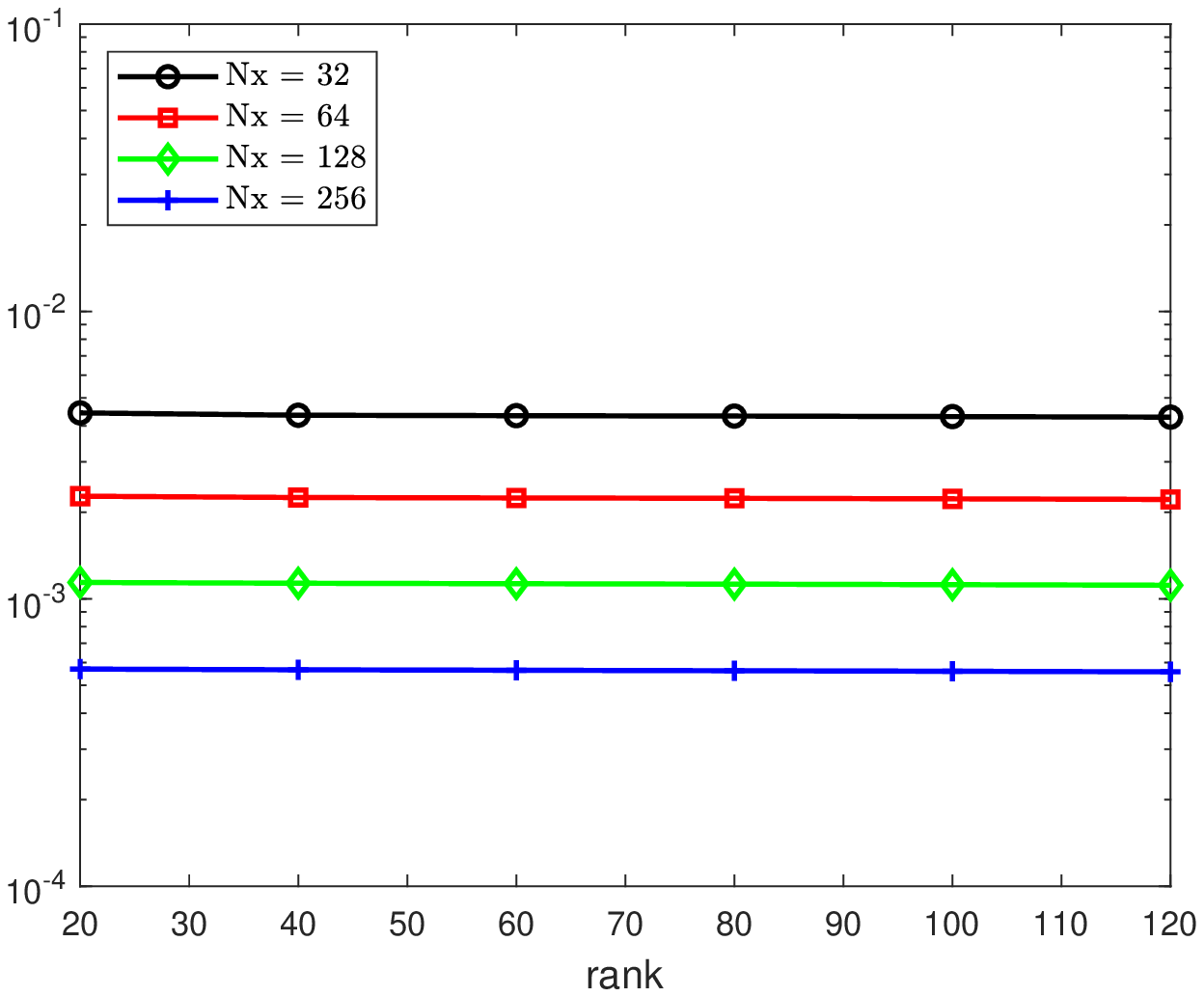}
\includegraphics[width=1.9in]{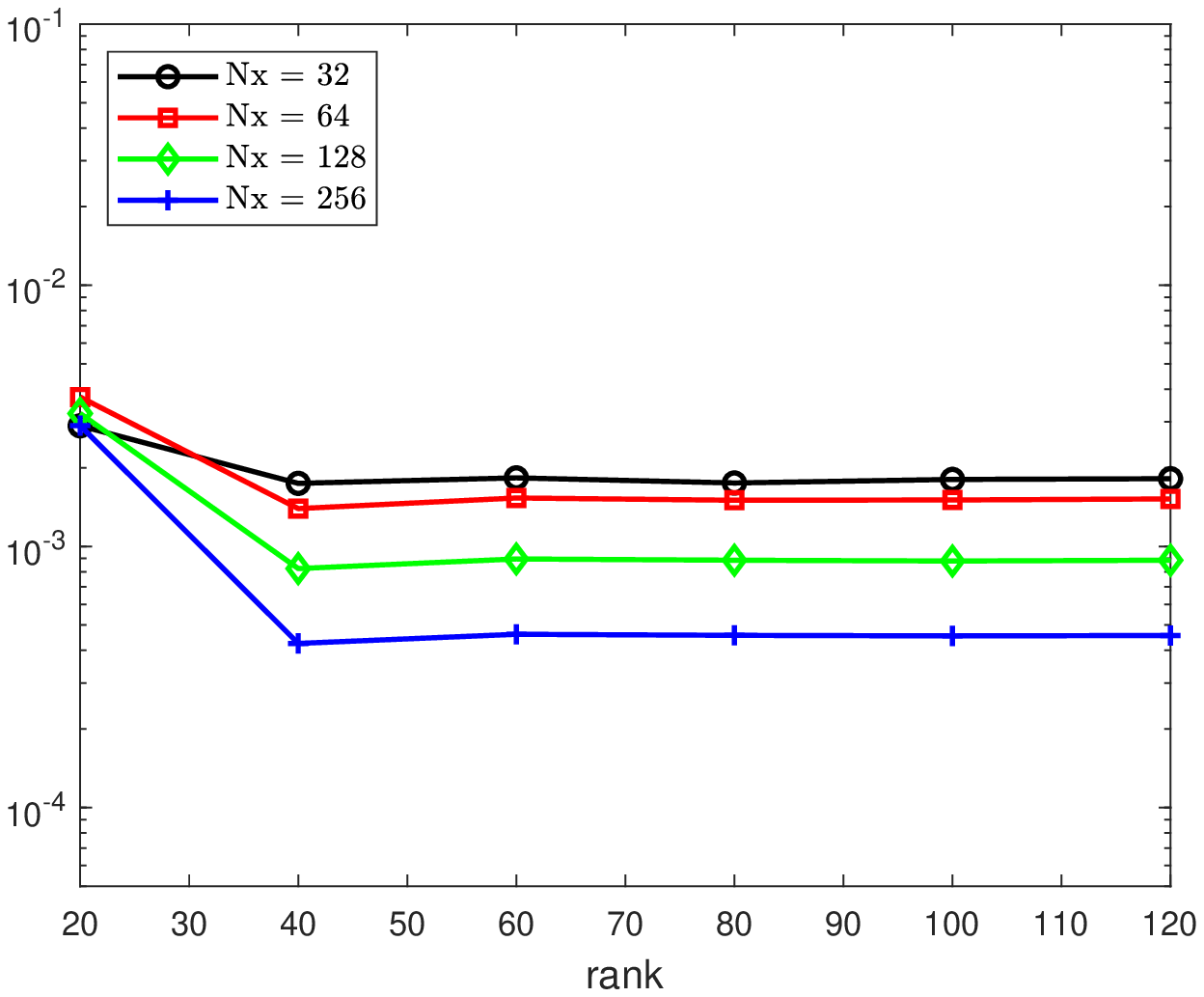}
\includegraphics[width=1.9in]{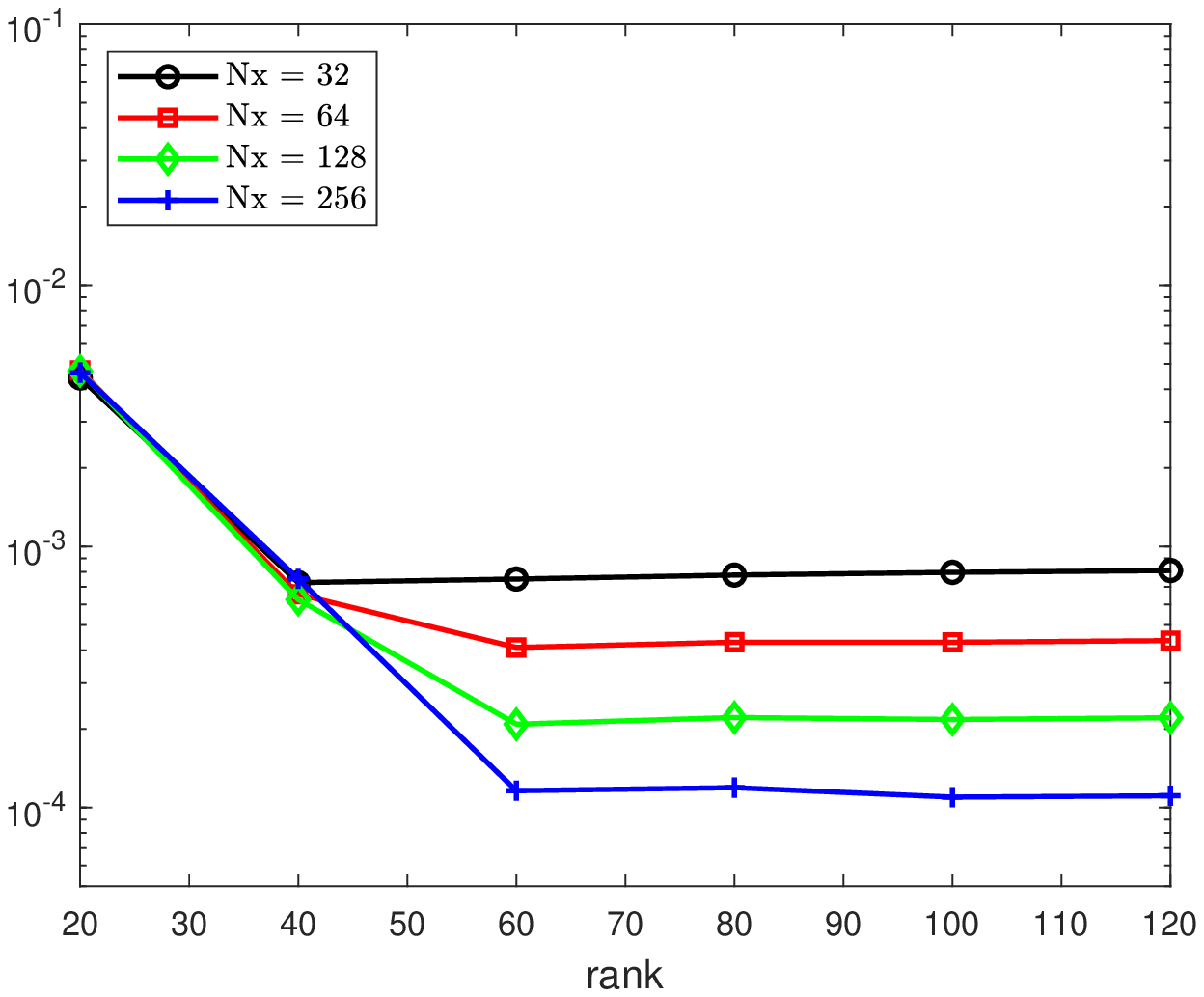}
\caption{Section~\ref{subsec:variable}: variable scattering coefficient. Density profile with $y = 0$ of the low rank solution and full tensor solution on a $256\times 256$ mesh at time $t=0.002$ (top left), $t=0.006$ (top middle), and $t=0.010$ (top right). Difference (\ref{diff}) between the low rank solution and full tensor solution computed on different meshes and with different ranks at time $t=0.002$ (bottom left), $t=0.006$ (bottom middle), and $t=0.010$ (bottom right).}
\label{figure:4}
\end{center}
\end{figure}

The comparison of the low rank solution and the full tensor solution on a $256\times 256$ mesh at different times is shown in Figure~\ref{figure:4} (top). We can see that the low rank solution matches well with the full tensor solution except for rank $r=20$. To quantitatively understand the rank dependence, we compute the difference of two solutions on the same mesh as follows
\begin{equation} \label{diff}
\left(\Delta x^2\sum_{k,l=1}^{N_x} \left(\rho_{\text{low rank}}(x_{k+\frac{1}{2}},y_{l+\frac{1}{2}})-\rho_{\text{full tensor}}(x_{k+\frac{1}{2}},y_{l+\frac{1}{2}})\right)^2\right)^{\frac{1}{2}}.
\end{equation}
and track how this evolves in time under certain fixed ranks $r$ ranging from $20$ to $120$. The results are shown in Figure~\ref{figure:4} (bottom). The common trend is that once the rank is increased to a certain level, the difference saturates. This is because then the spatial error dominants. Also it is clear that the rank of the solution in this problem increases gradually with time.

\begin{figure}[!htp]
\begin{center}
\includegraphics[width=2.6in]{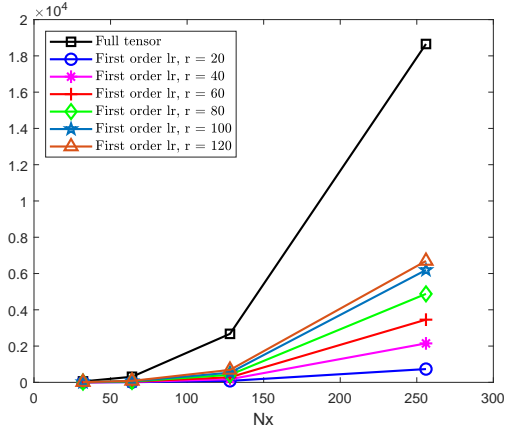}
\caption{Section~\ref{subsec:variable}: variable scattering coefficient. Computational time (in seconds) needed for the low rank method and full tensor method to compute the solution at time $t = 0.012$.}
\label{figure:5}
\end{center}
\end{figure}

In addition, we record the computational time needed to compute the solution to $t = 0.012$ for both methods on an i7-8700k @3.70 GHz CPU in Figure~\ref{figure:5}. The speedup of the low rank method is significant, especially for a large number of spatial points $N_x$.

\subsection{Two-material test}
\label{subsec:two}

\begin{figure}[!htp]
\begin{center}
\includegraphics[width=2.6in]{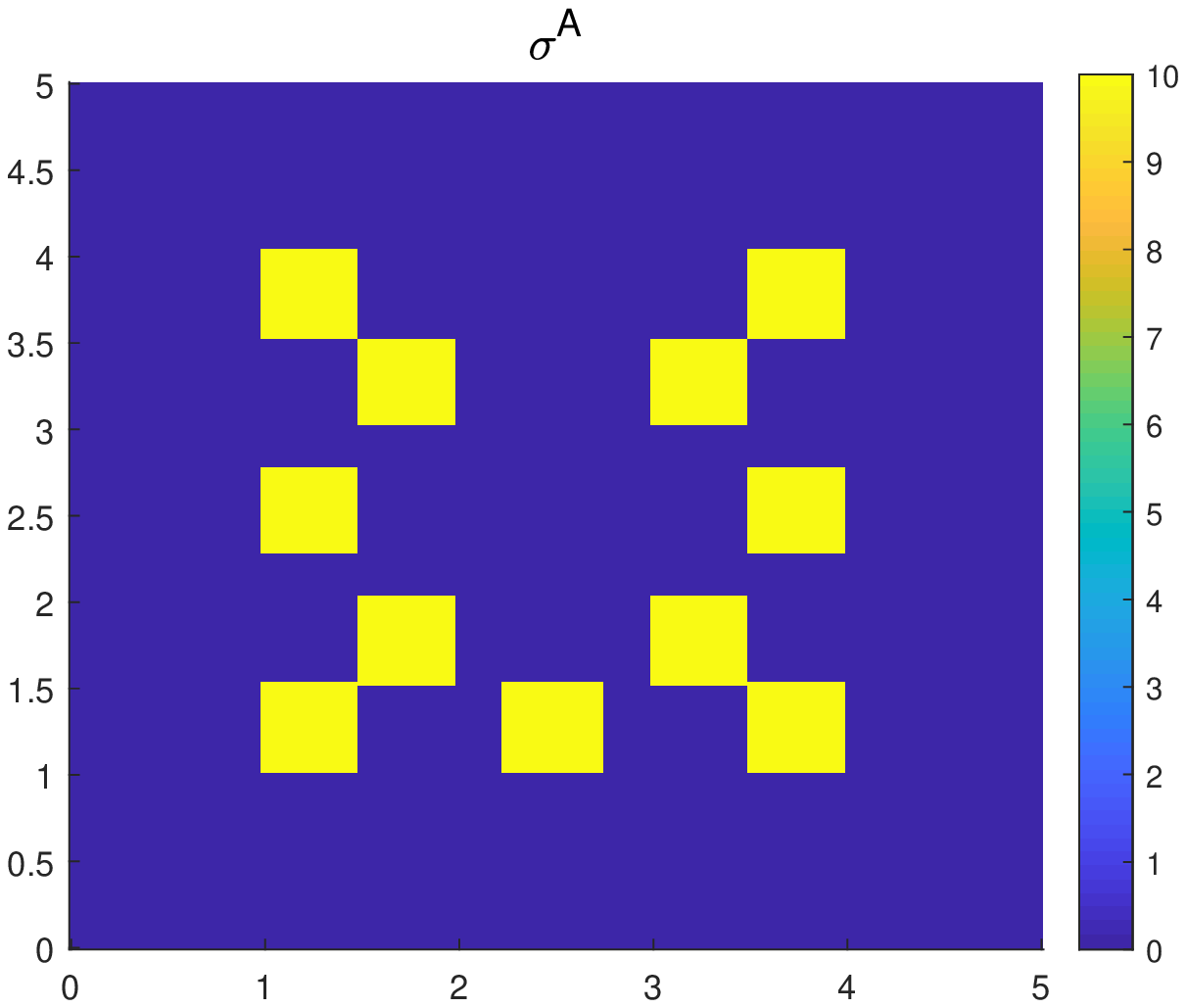}
\includegraphics[width=2.6in]{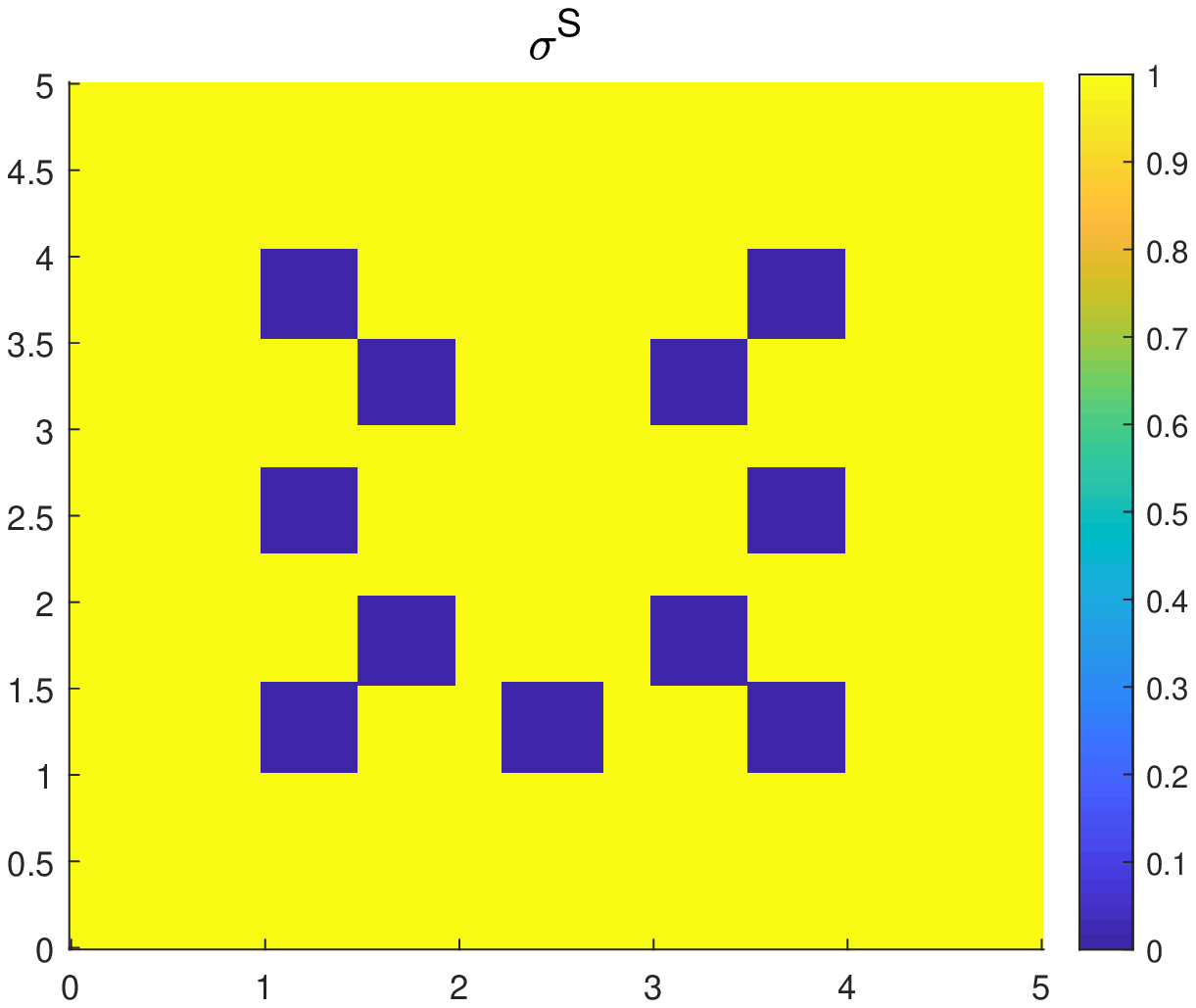}
\caption{Section~\ref{subsec:two}: two-material test. Profiles of absorption coefficient $\sigma^A$ (left) and scattering coefficient $\sigma^S$ (right). Each square block in the computational domain is a $0.5 \times 0.5$ square. In the left figure, yellow square blocks represent that $\sigma^A = 10$ and for the rest blue region $\sigma^A = 0$; in the right figure, blue square blocks represent that $\sigma^S = 0$ and for the rest yellow region $\sigma^S = 1$.}
\label{figure:7}
\end{center}
\end{figure}

The two-material test models a domain with different materials with discontinuities in material cross sections and source term. It is a slight modification of the lattice benchmark problem for linear transport equation. Here we choose the computational domain as $[0,5]^2$ with the absorption coefficient $\sigma^A$ and scattering coefficient $\sigma^S$ given as in Figure~\ref{figure:7}. The source term is given by
\be
	G(x,y) = 
	\left\{
	\bal
	&1,& &(x,y) \in [2,3]^2,\\
	&0,& &\text{otherwise.}
	\eal
	\right.
\ee

\begin{figure}[!htp]
\begin{center}
\includegraphics[width=2.5in]{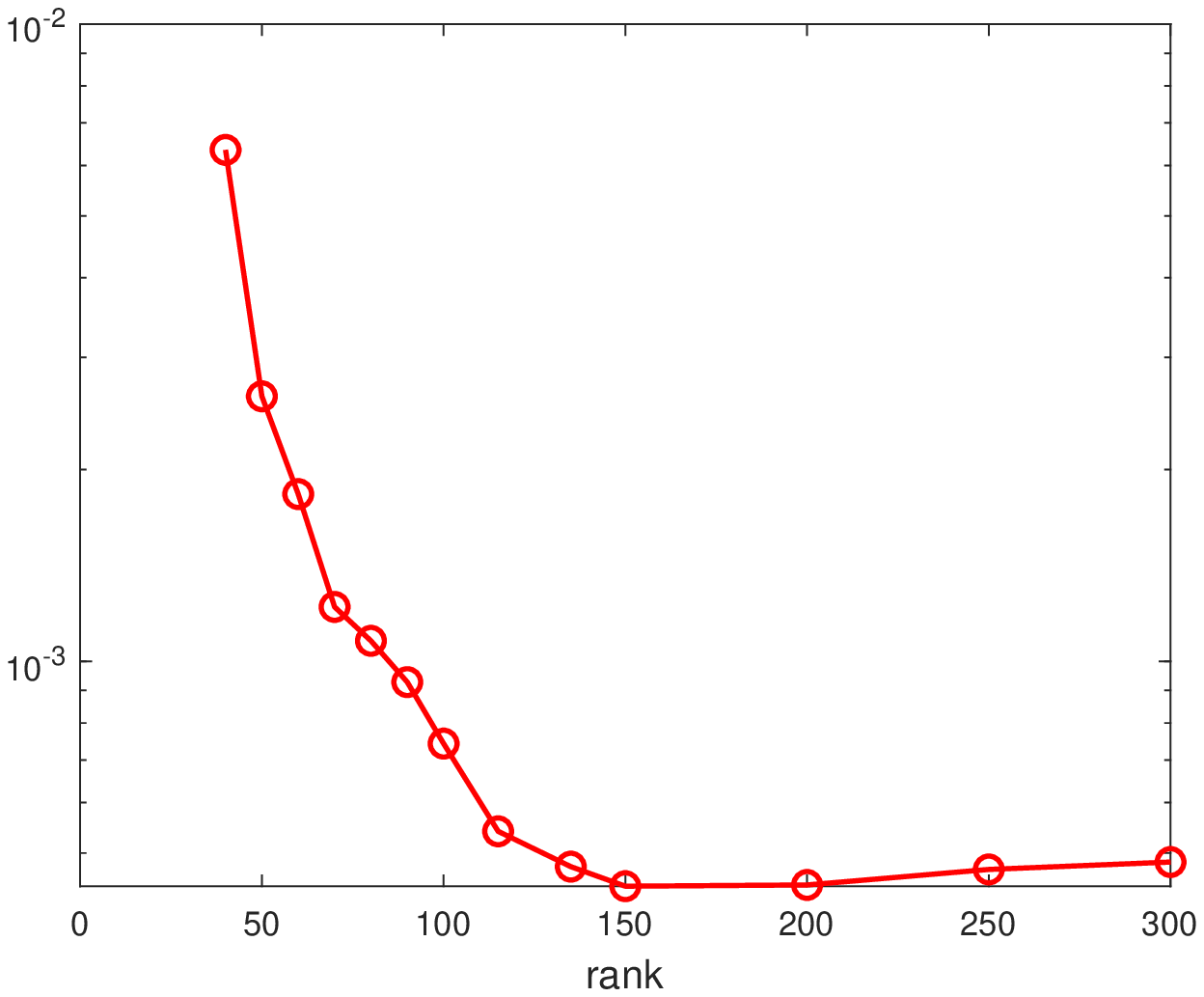}
\includegraphics[width=2.7in]{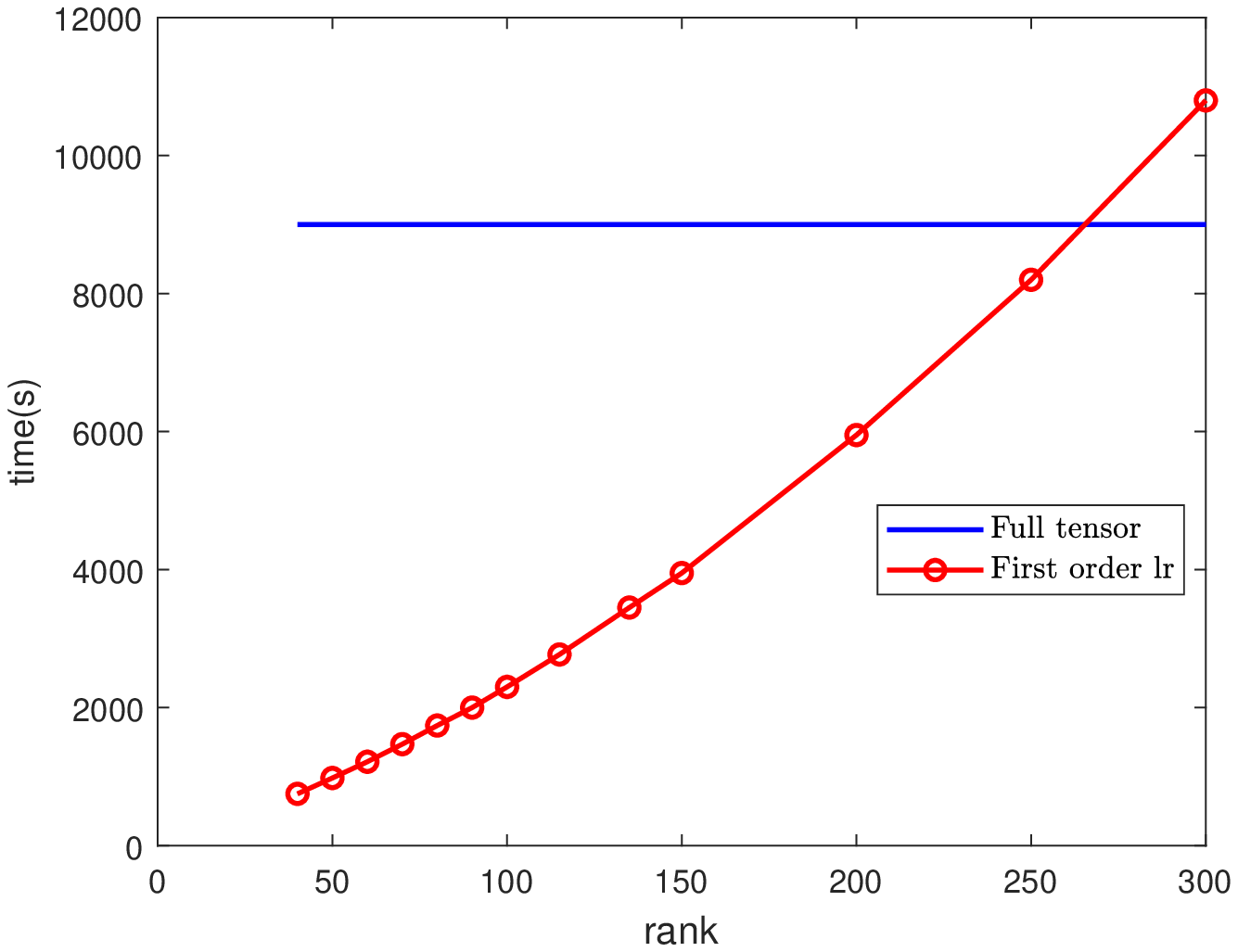}
\caption{Section~\ref{subsec:two}: two-material test ($\varepsilon = 1$). Difference (\ref{diff}) between the low rank solution with different ranks and full tensor solution at time $t=1.7$ (left). Computational time (in seconds) needed for the low rank method with different ranks and full tensor method to compute the solution at $t=1.7$ (right).}
\label{figure:9}
\end{center}
\end{figure}

\begin{figure}[!htp]{}
\begin{center}
\includegraphics[width=2.6in]{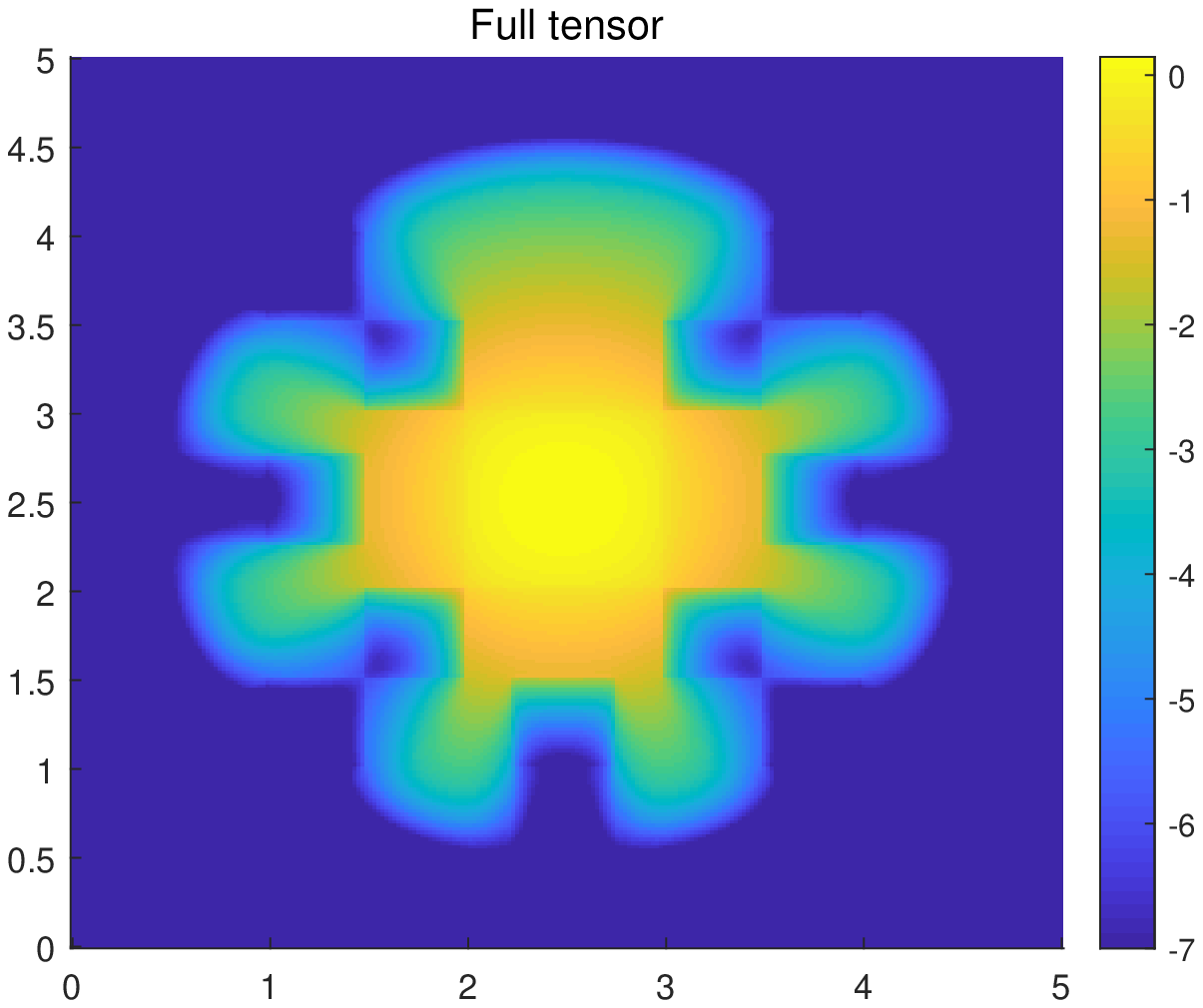}
\includegraphics[width=2.6in]{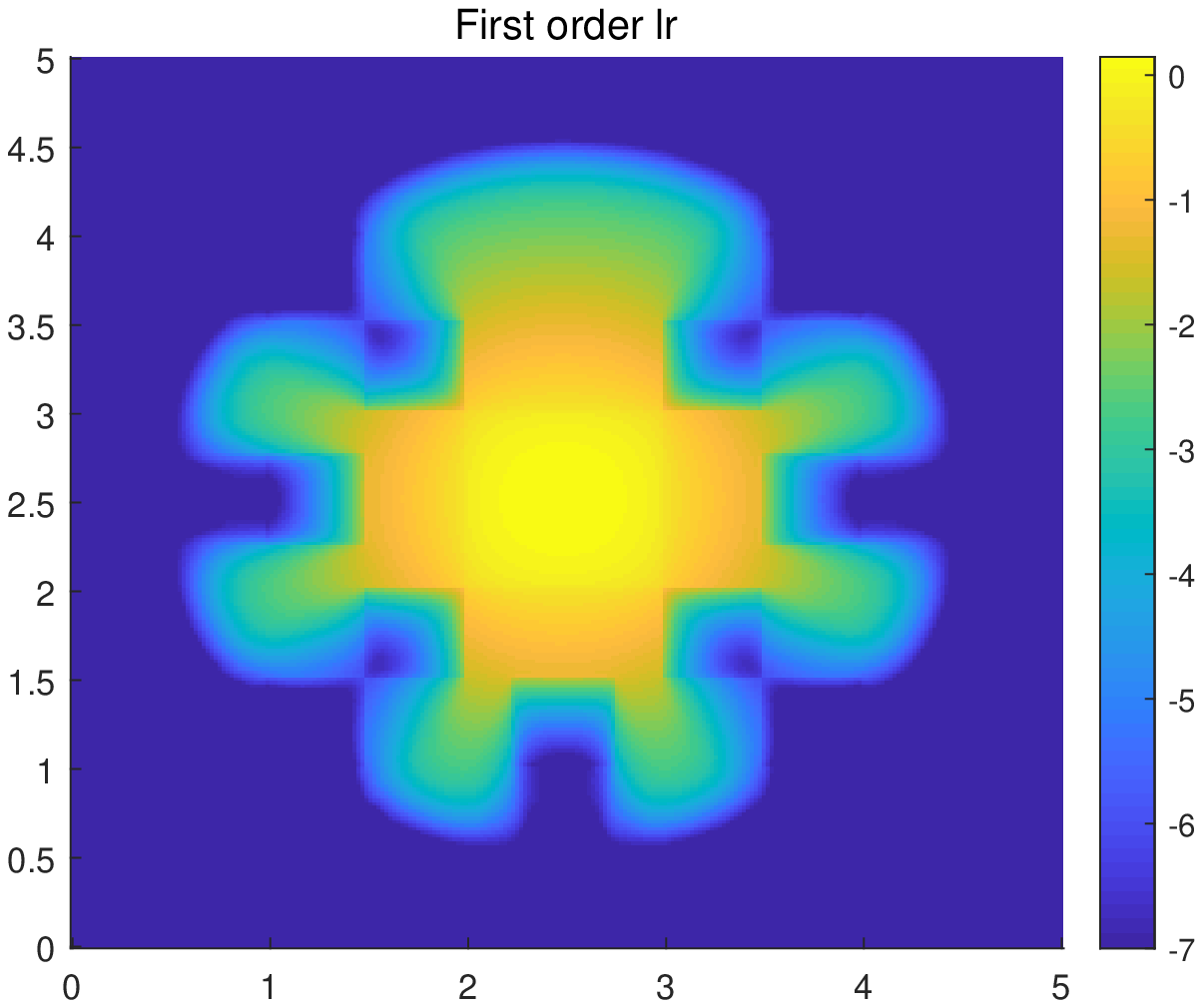}
\includegraphics[width=2.6in]{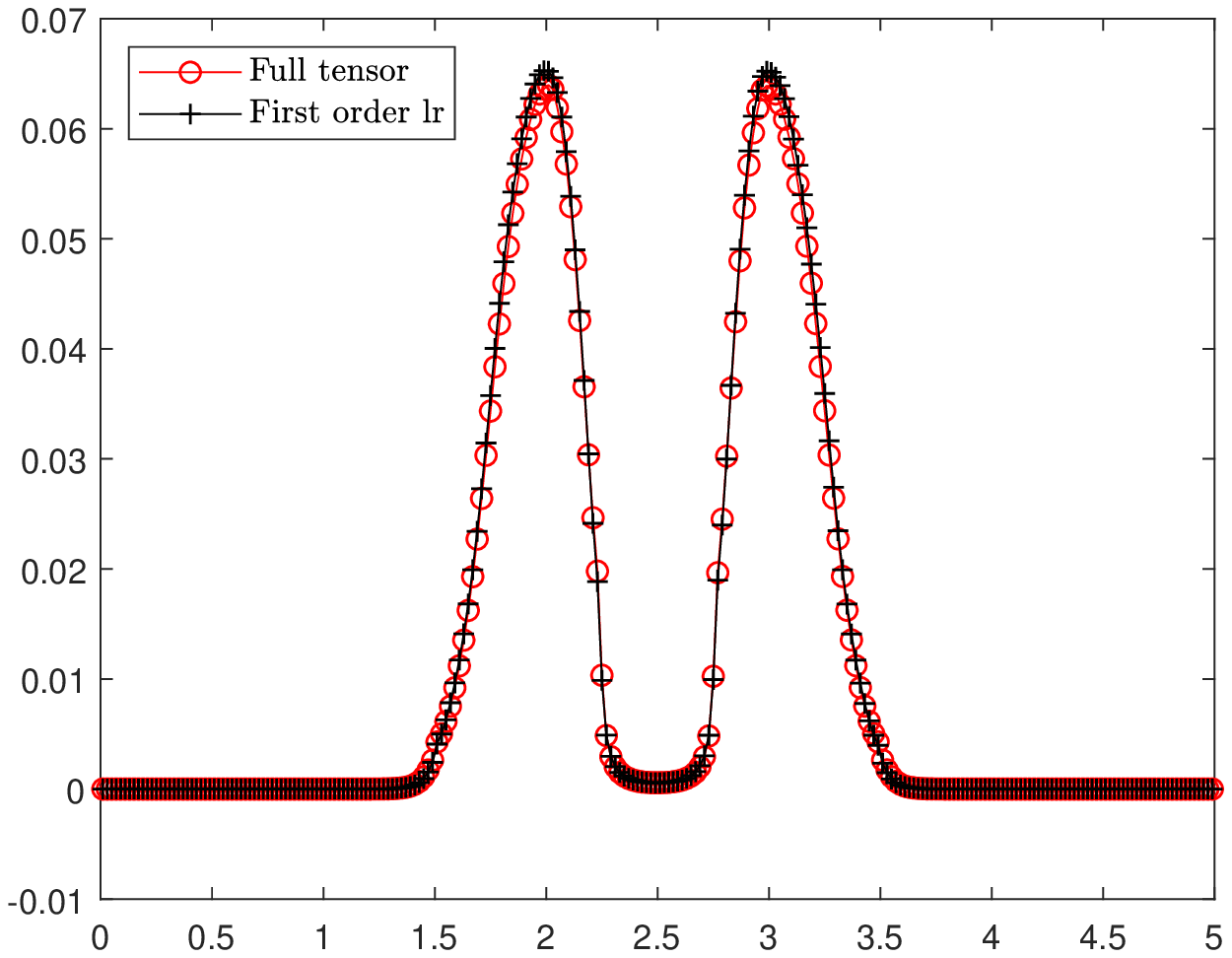}
\includegraphics[width=2.6in]{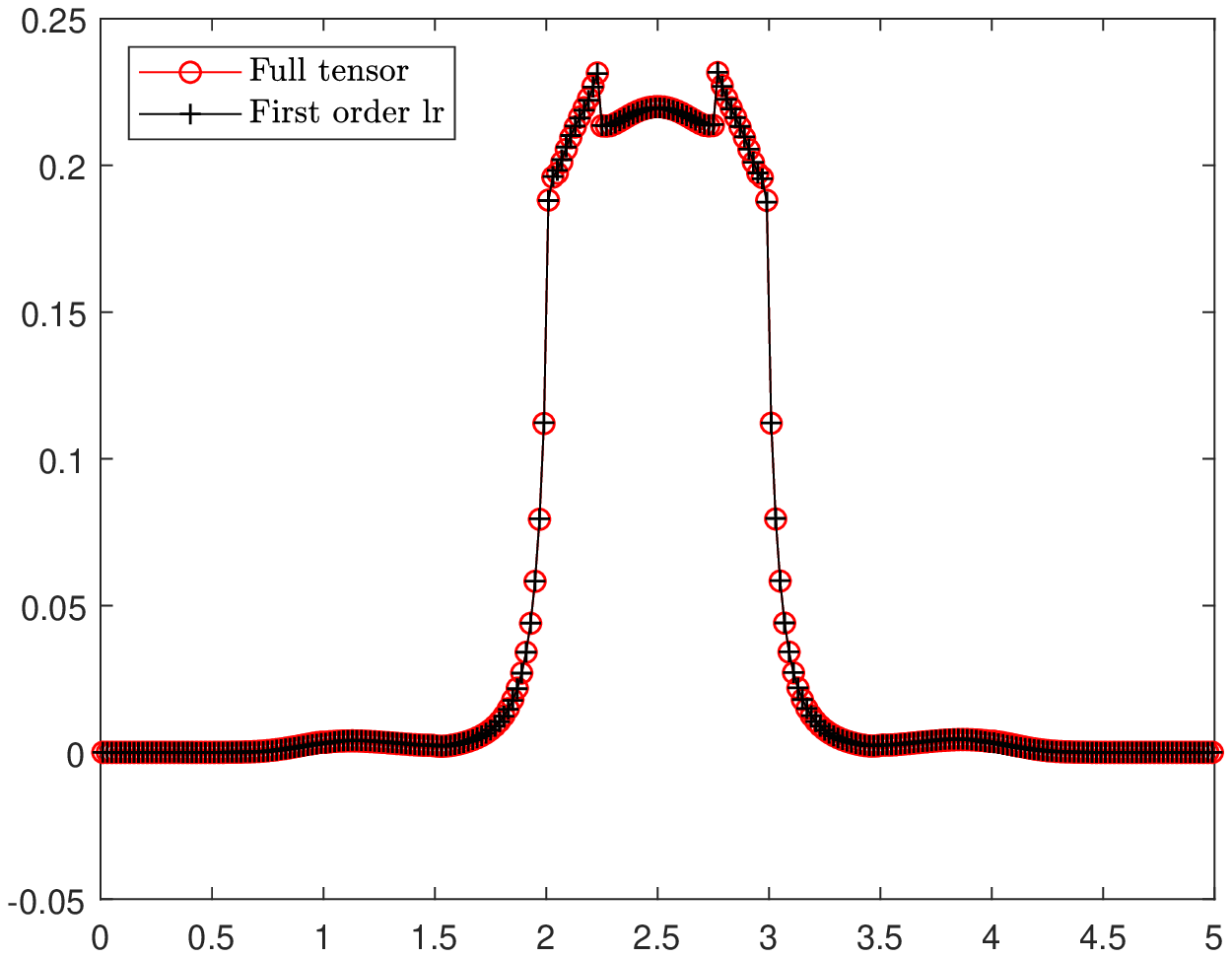}
\includegraphics[width=2.6in]{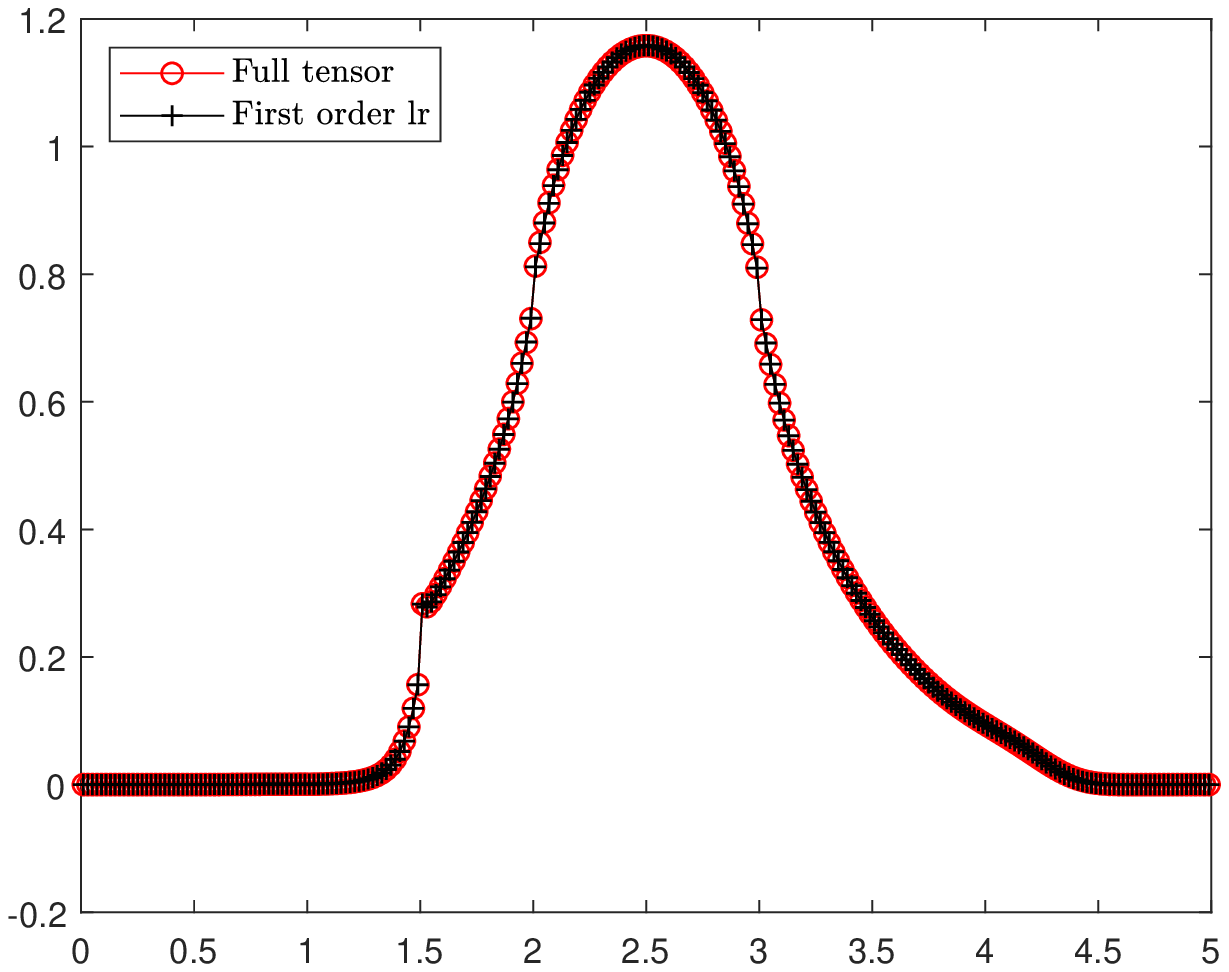}
\includegraphics[width=2.6in]{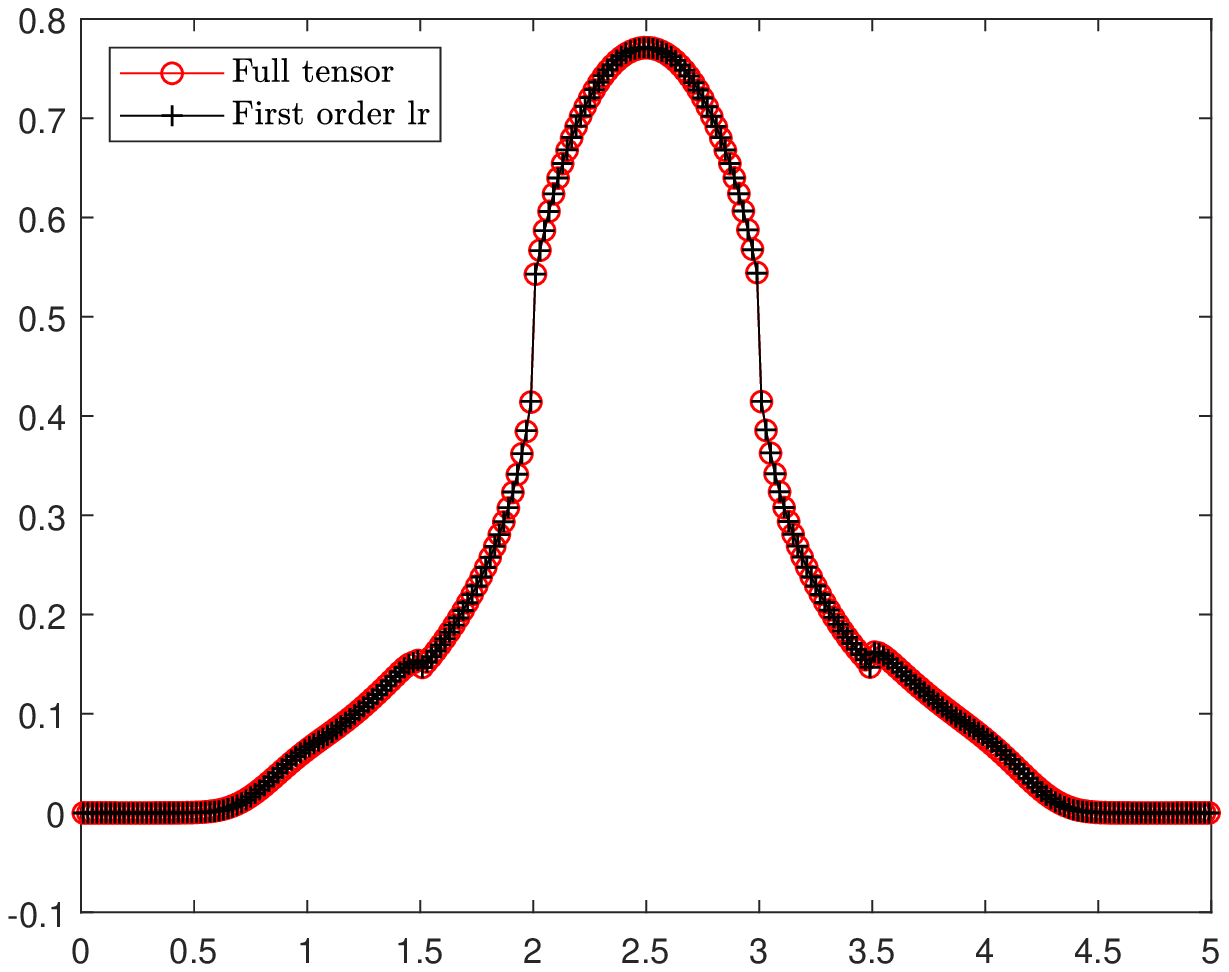}
\caption{Section~\ref{subsec:two}: two-material test ($\varepsilon = 1$). Contour plot of the log density at time $t=1.7$ of the full tensor solution (top left) and low rank solution (top right) on a $250\times250$ mesh. Density slice of both solutions along $x = 1$ (middle left), $x = 1.5$ (middle right), $x = 2.5$ (bottom left), and $x = 3$ (bottom right). $r = 150$ in the low rank method.}
\label{figure:8}
\end{center}
\end{figure}

We set $\varepsilon = 1$ and compare the first order low rank method with the first order full tensor method. For both methods, we choose $N_x=N_y=250$, $N_{\vv} = 2702$ Lebedev quadrature points on $\mathbb{S}^2$, and same mixed CFL condition $\dt t= 0.1\min(\sigma^S)\dt x^2 + 0.1\varepsilon\dt x $. The initial condition is given by
\be
f(t=0,x,y,\xi,\eta,\gamma) = \fl {1}{4\pi\varsigma^2}\exp\left(-\fl {(x-2.5)^2+(y-2.5)^2}{4\varsigma^2}\right), \quad \varsigma^2 = 10^{-2},\quad  (x,y) \in [0,5]^2.
\ee
We test different ranks from $40$ to $300$ in the low rank method and compare it with the full tensor solution. The error and computational time are reported in Figure~\ref{figure:9}. It is clear that at around rank $r = 150$, the spatial error dominates and increasing the rank further will have no gain in solution accuracy. Moreover, at $r=150$, the efficiency of the low rank method is clearly better than the full tensor method. We then fix $r=150$ and plot both the low rank solution and full tensor solution at $t=1.7$ in Figure~\ref{figure:8}, where a good match is obtained.

\begin{figure}[!htp]
\begin{center}
\includegraphics[width=2.6in]{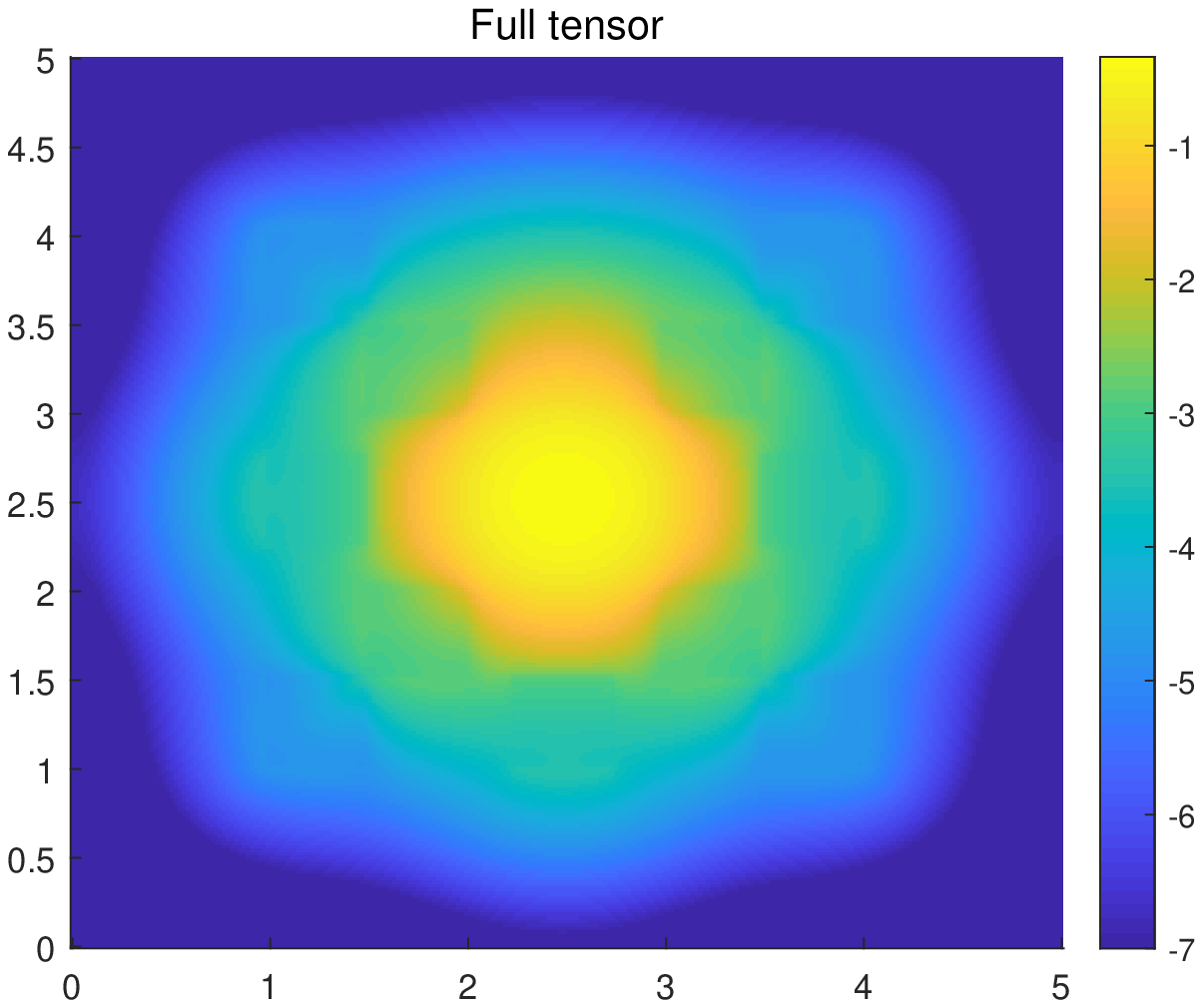}
\includegraphics[width=2.6in]{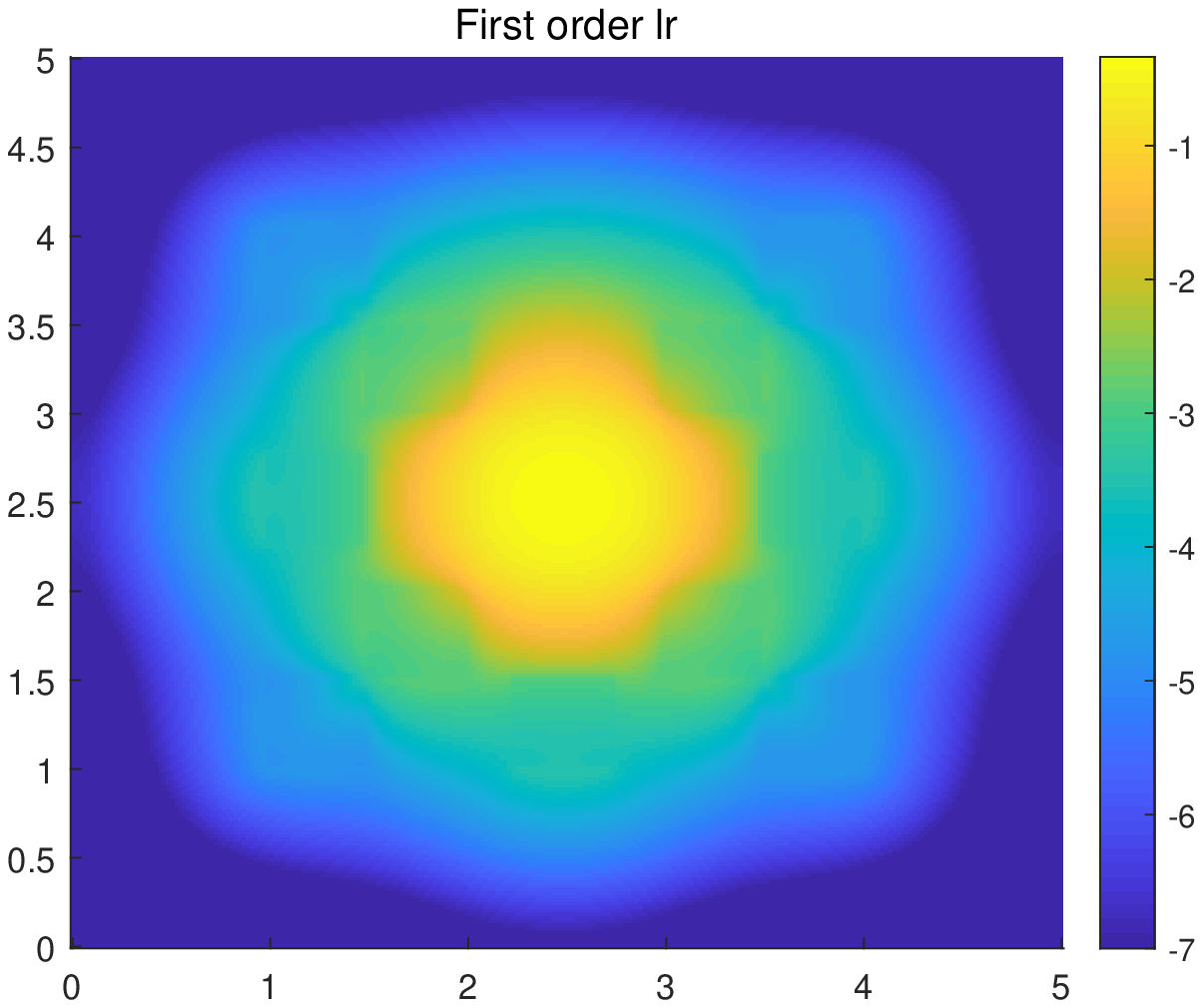}
\includegraphics[width=2.6in]{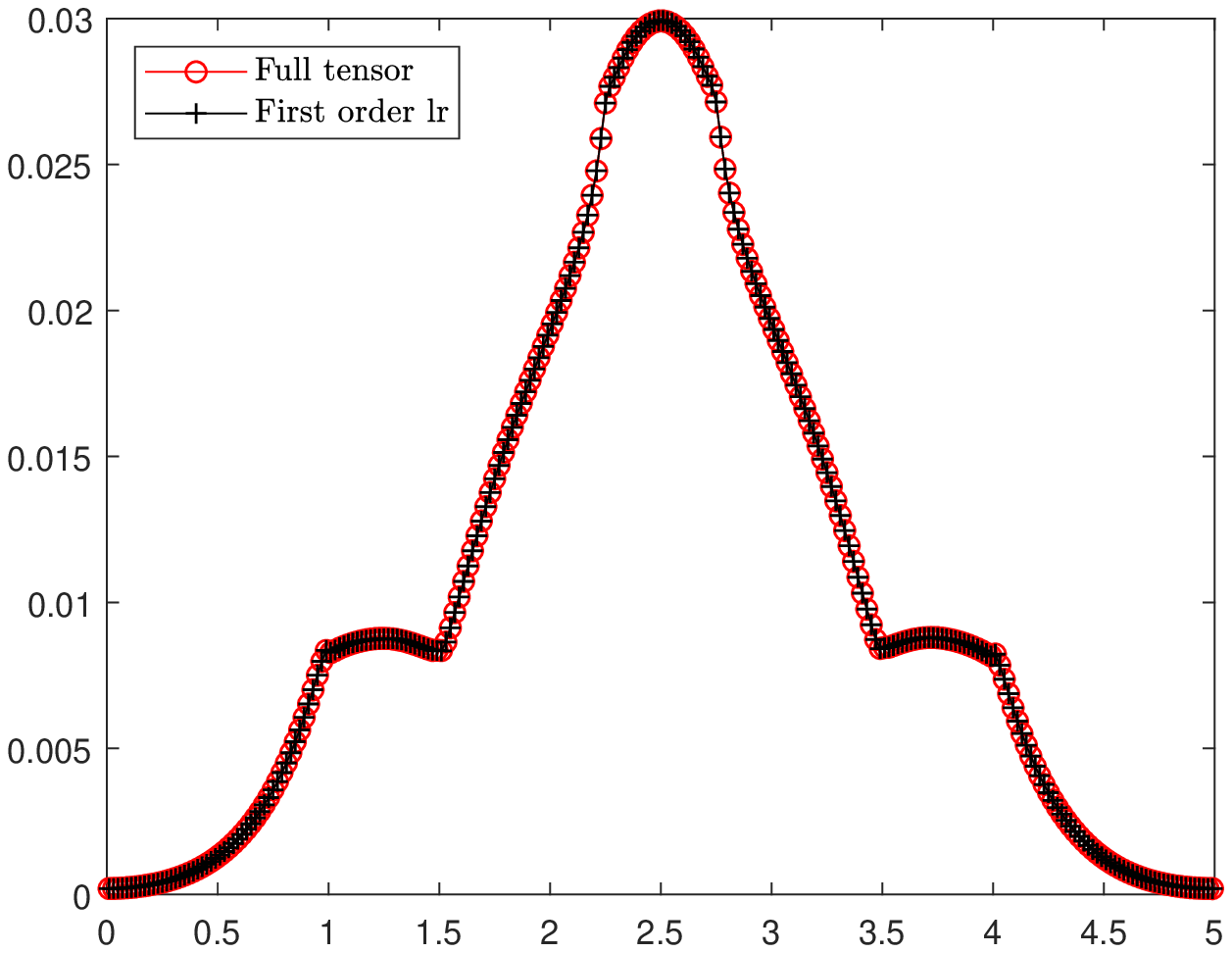}
\includegraphics[width=2.6in]{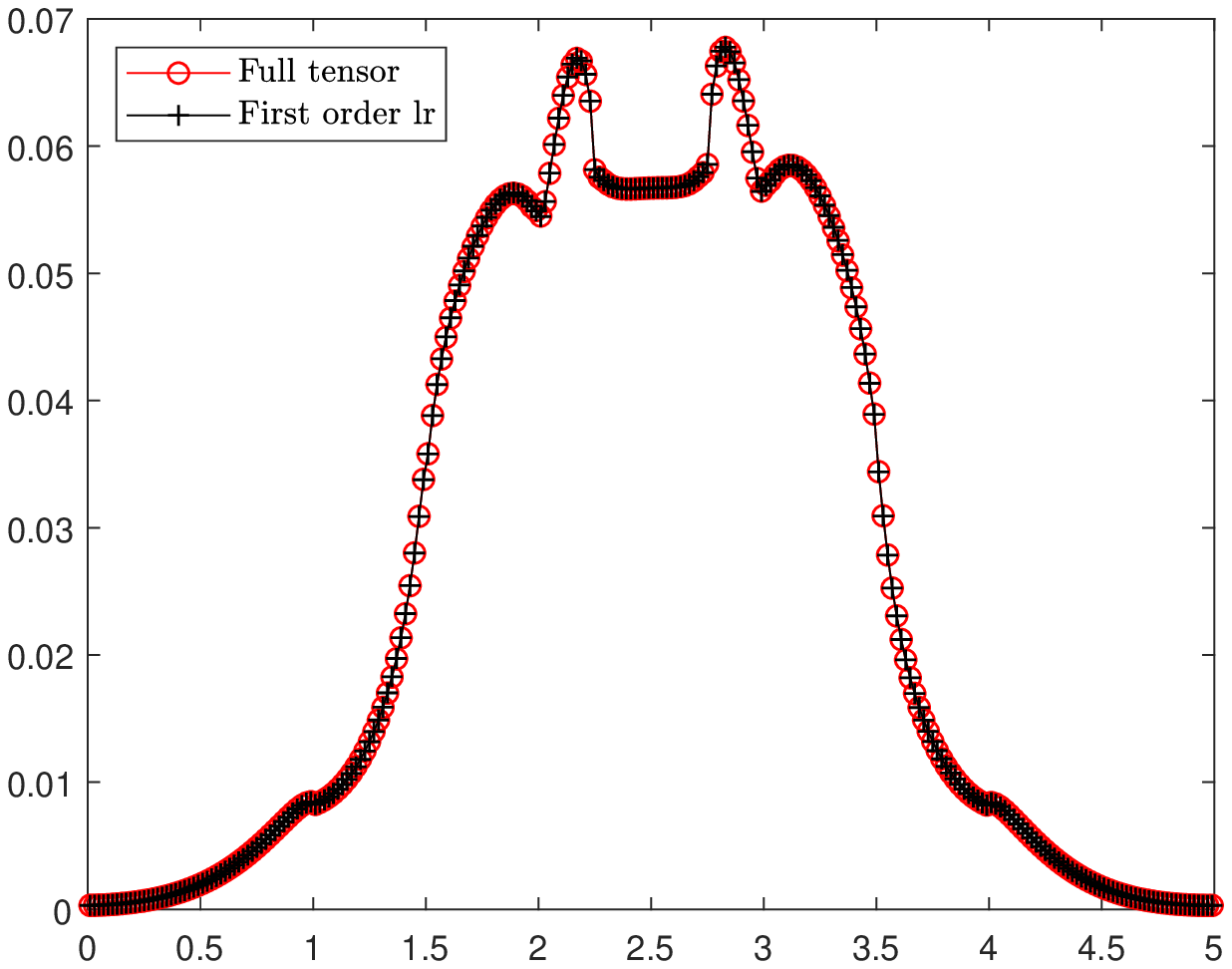}
\includegraphics[width=2.6in]{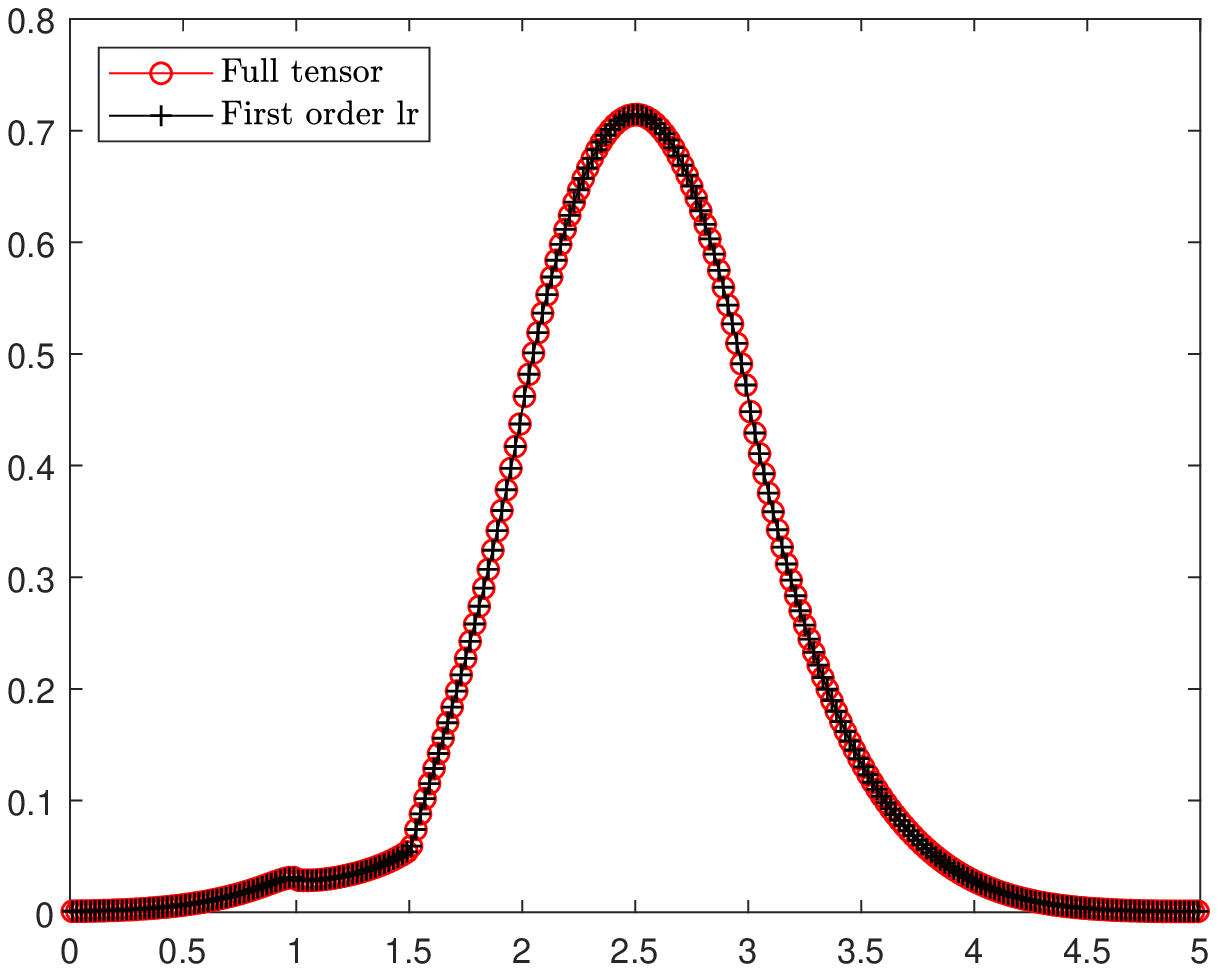}
\includegraphics[width=2.6in]{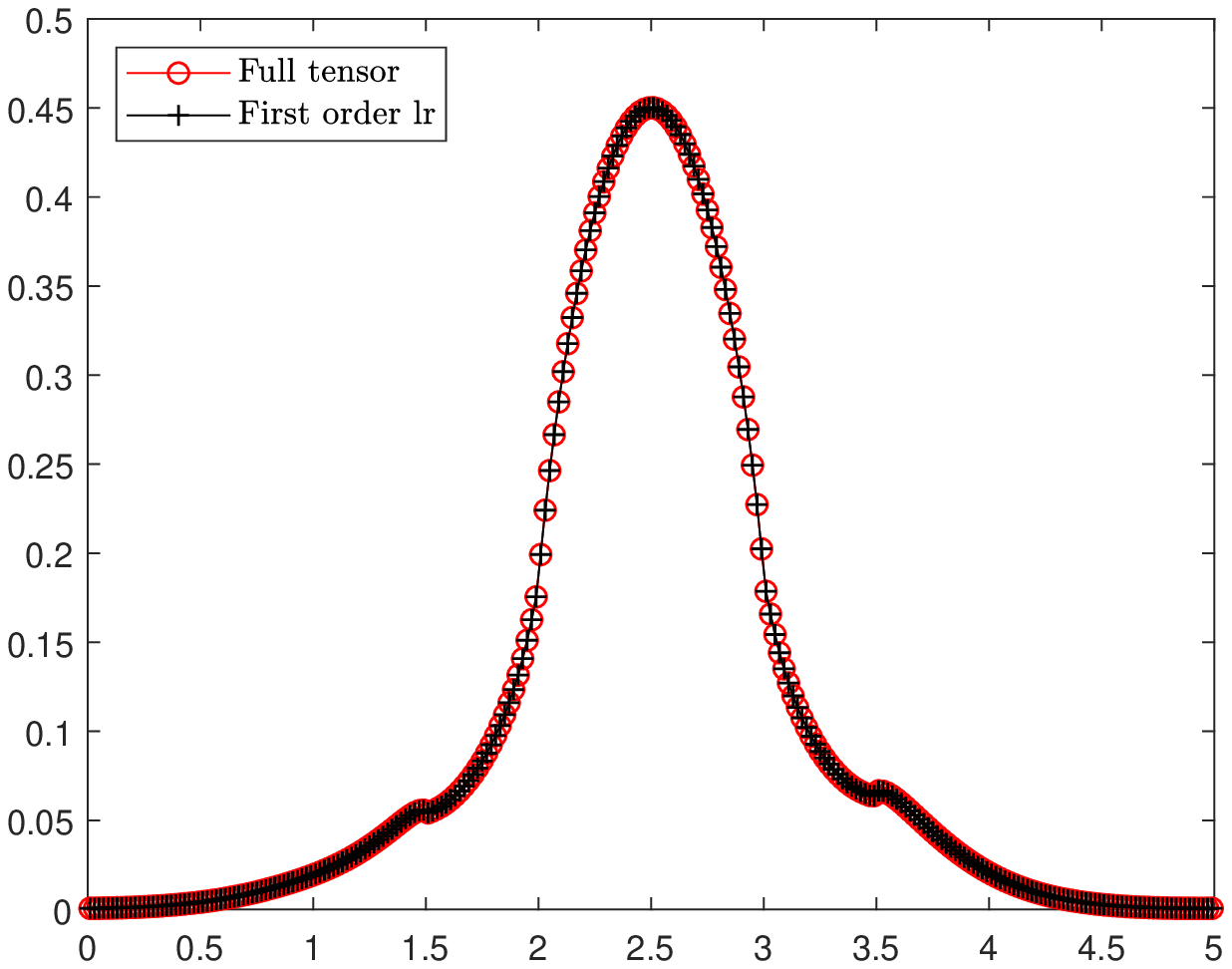}
\caption{Section~\ref{subsec:two}: two-material test ($\varepsilon = 0.1$). Contour plot of the log density at time $t=0.6$ of the full tensor solution (top left) and low rank solution (top right) on a $250\times250$ mesh. Density slice of both solutions along $x = 1$ (middle left), $x = 1.5$ (middle right), $x = 2.5$ (bottom left), and $x = 3$ (bottom right). $r = 100$ in the low rank method.}
\label{figure:10}
\end{center}
\end{figure}

In addition, we consider another scenario with $\varepsilon = 0.1$. The same parameters are used as in the case of $\varepsilon=1$, except we set the rank $r = 100$ in the low rank method (because we expect the rank of the solution to decrease as $\varepsilon$ decreases). The solutions of the low rank method and full tensor method at time $t=0.6$ are shown in Figure~\ref{figure:10}, where we again observe good agreement. An optimal (and possibly smaller) rank can be determined similarly as in Figure~\ref{figure:9}, we omit the result.

\subsection{Line source test}
\label{subsec:line}

\begin{figure}[!htp]
\begin{center}
\includegraphics[width=2in]{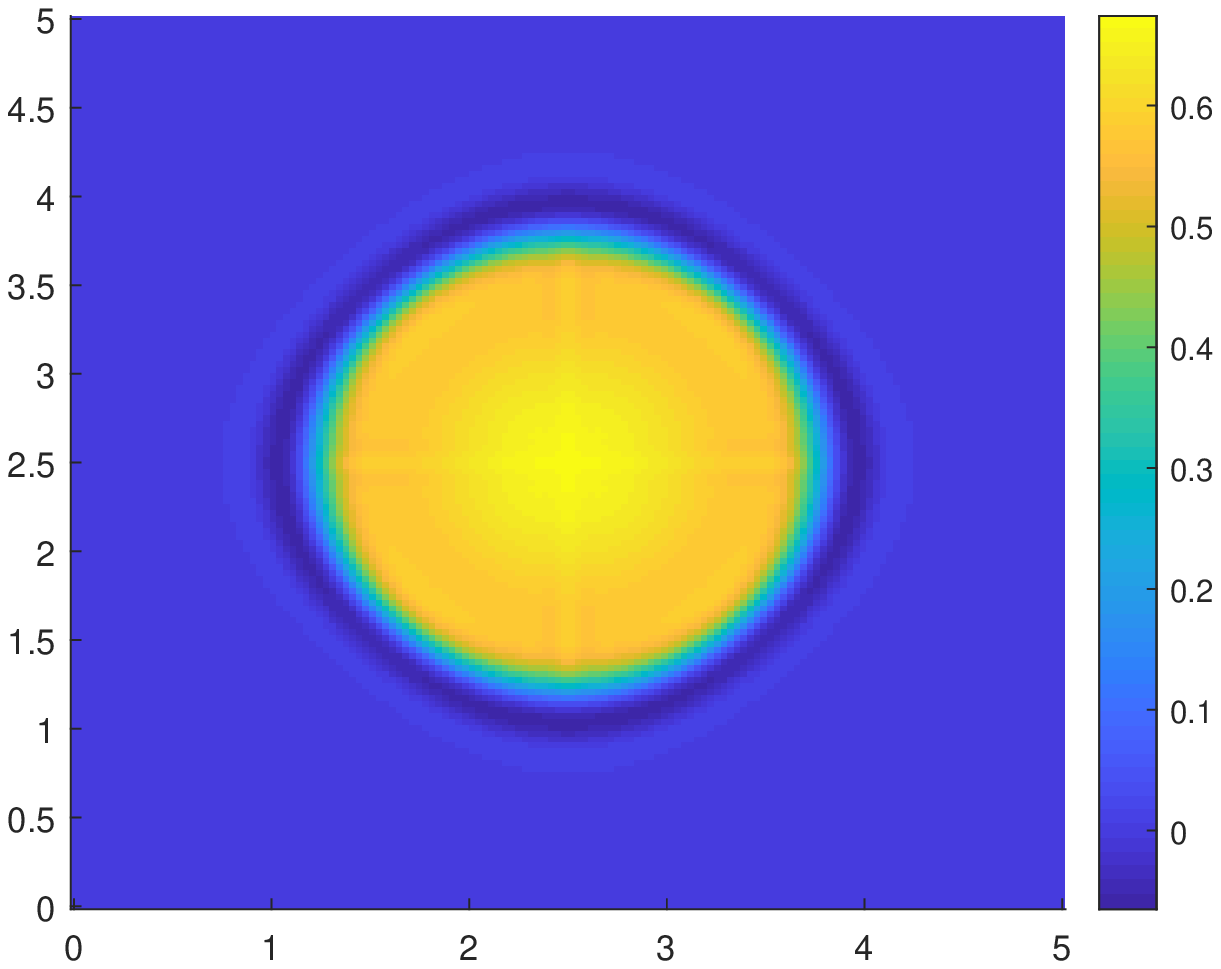}
\includegraphics[width=2in]{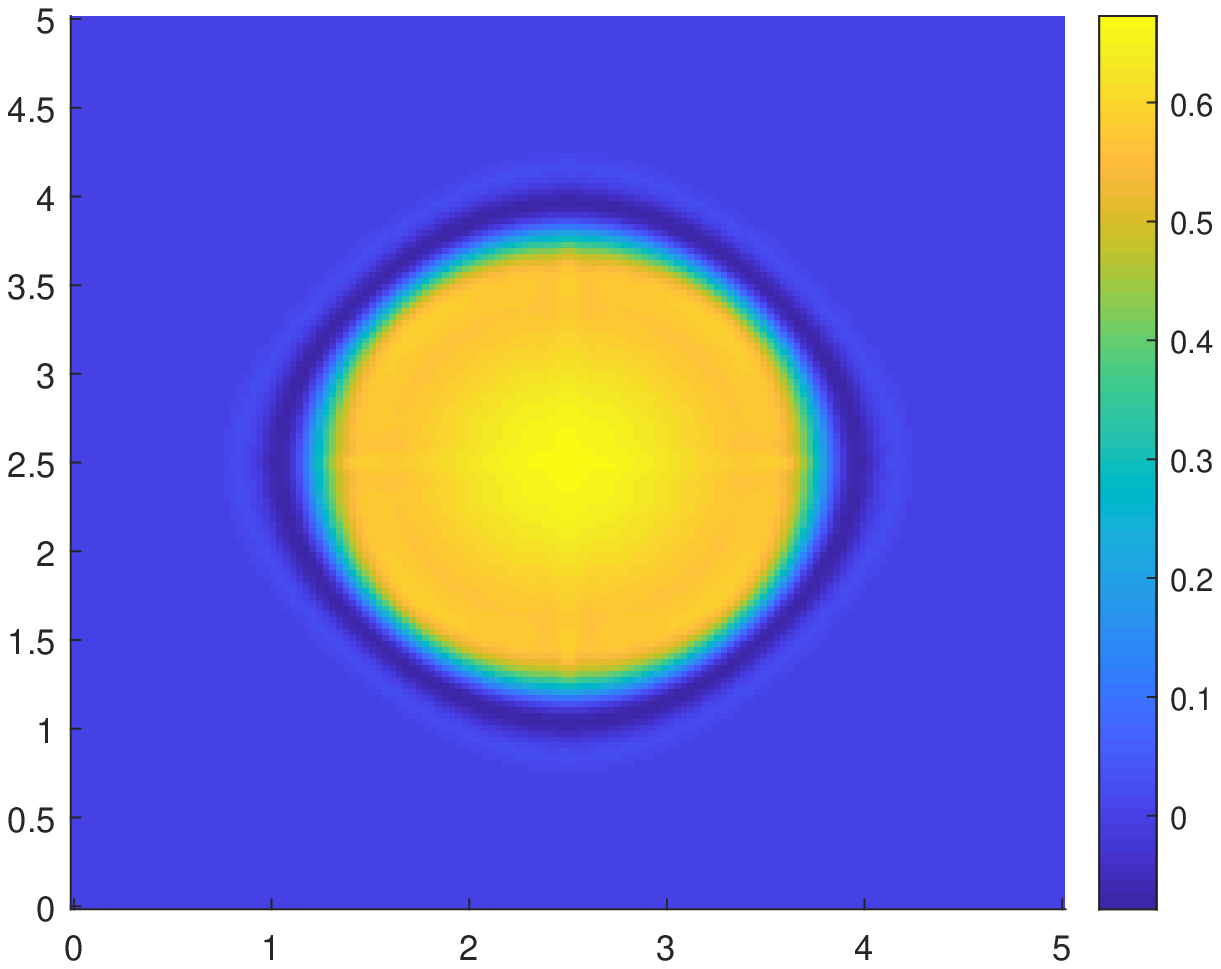}
\includegraphics[width=2in]{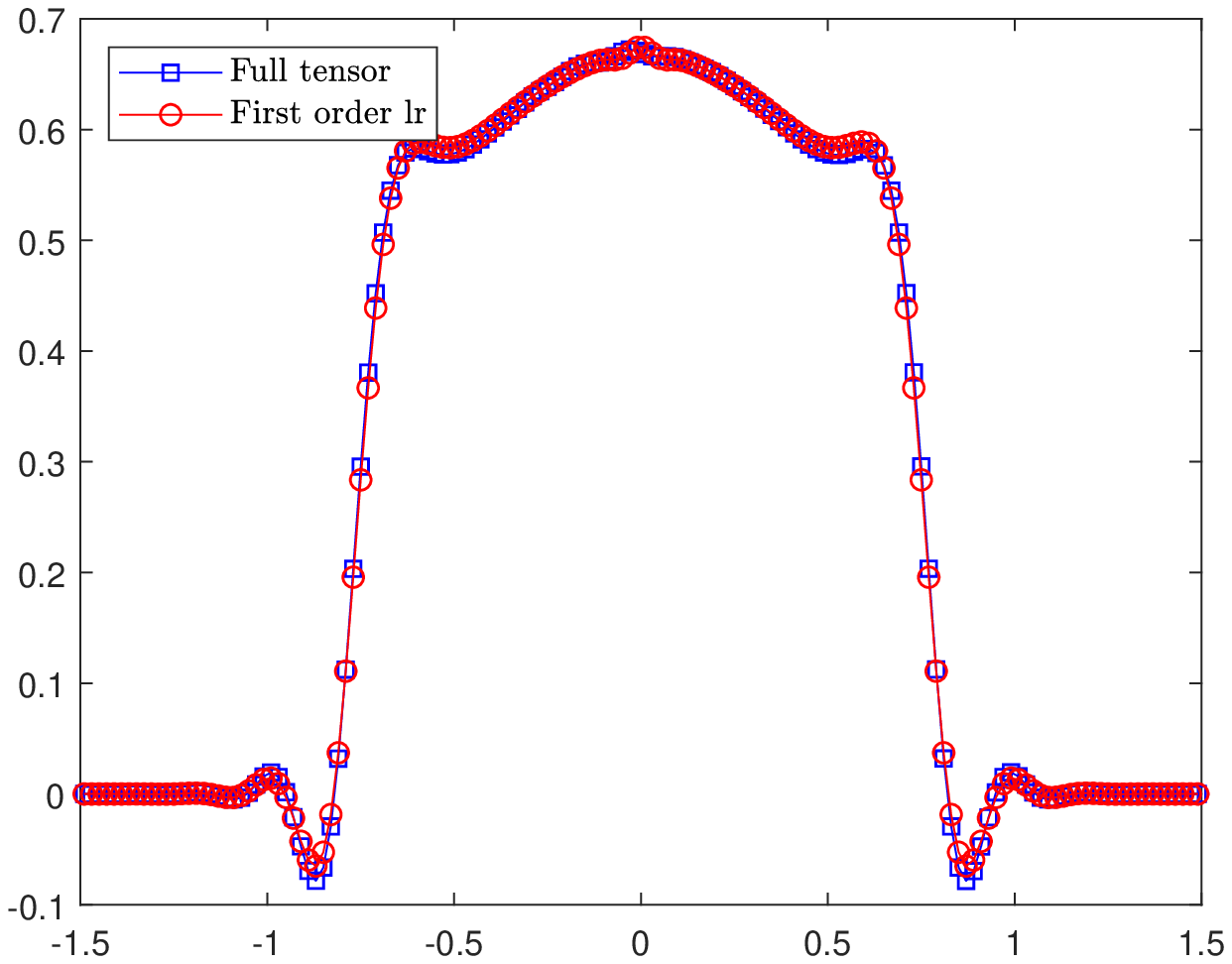}
\caption{Section~\ref{subsec:line}: line source test. Density profile of the full tensor solution (left) and low rank solution (middle) on a $150 \times 150$ mesh, and comparison of two solutions along $y = 0$ (right) at time $t=0.7$. $r = 600$ in the low rank method.}
\label{figure:6}
\end{center}
\end{figure}

We finally consider the line source test which is another important benchmark test for the linear transport equation. Here we approximate the initial delta function via (\ref{eqn:initial_gauss}) with a much smaller $\varsigma^2 = 4\times10^{-4}$. $\sigma^S=1$ and $\sigma^A=G=0$. We set $\varepsilon = 1$ and compare the first order low rank method with the full tensor method. For both methods, we choose the computational domain as $[-1.5,1.5]^2$ with $N_x=N_y=150$, $N_{\vv} = 5810$ Lebedev quadrature points on $\mathbb{S}^2$, and the same mixed CFL condition $\dt t= 0.025\dt x^2 + 0.025\varepsilon\dt x $. We fix the rank as $r = 600$ in the low rank method. The density profiles of both methods at time $t = 0.7$ are shown in Figure \ref{figure:6}. We can see that the solutions match well. 

We would like to mention that this is a difficult problem compared to the cases considered previously. Many more points are need on the sphere to get a reasonable solution. Nevertheless, there are still oscillations in the solution (for both the full tensor and the low rank method). This is a well-known artifact in the $S_N$ method. In addition, we found that a higher rank and a more stringent CFL condition is needed in the low rank method. We believe part of the reason are the numerical oscillations, which can be tempered by applying a proper filter or using a positivity-preserving scheme. We refer to \cite{LFH19}, and references therein, for more details.

\section{Conclusion}
\label{sec:con}

We have introduced a dynamical low-rank method for the multi-scale multi-dimensional linear transport equation. The method is based on a macro-micro decomposition of the equation and uses the low rank approximation only for the micro part of the solution. The key feature of the proposed scheme is that it is explicitly implementable, asymptotic-preserving in the diffusion limit, and maintains second order in both kinetic and diffusive regimes. A series of numerical examples in 2D including some well-known benchmark tests have been performed to validate the accuracy, efficiency, rank dependence, and AP property of the proposed method. Some interesting ongoing and future work includes adaptive rank selection and the theoretical investigation of rank dependence of the solution in the kinetic regime.


\bibliographystyle{plain}
\bibliography{hu_bibtex,le_bibtex}


\end{document}